\documentclass[11pt]{article}

\newcommand{\documentdate}{18 April 2019}

\usepackage{a4wide,latexsym,varioref,amsmath}

\newcommand{\numsection}[1]{\section{#1}\setcounter{equation}{0}}

\renewcommand{\thefootnote}{(\arabic{footnote})}
\newcommand{\beqn}[1]{\begin{equation}\label{#1}}
\newcommand{\eeqn}{\end{equation}}
\newcommand{\req}[1]{(\ref{#1})}

\newcommand{\ms}{\;\;\;\;}
\newcommand{\tim}[1]{\;\; \mbox{#1} \;\;}
\newtheorem{theorem}{Theorem}[section]
\newtheorem{lemma}[theorem]{Lemma}

\newcommand{\llem}[2]{\vspace{\baselineskip} 
\noindent\framebox[\textwidth]{\parbox{0.95\textwidth}{
\begin{lemma} \label{#1} \rm #2 \end{lemma} } } \vspace{\baselineskip} }
\newcommand{\lthm}[2]{\vspace{\baselineskip} 
\noindent\framebox[\textwidth]{\parbox{0.95\textwidth}{
\begin{theorem} \label{#1} \rm #2 \end{theorem} } } \vspace{\baselineskip} }
\newlength{\thmw}
\newcommand{\lbthm}[3]{\vspace{\baselineskip}\noindent\hbox{%
  \lower\fboxrule\hbox{\vbox{\hrule\hbox{\vrule \kern-\fboxrule \vbox{%
  \vspace{\fboxsep} \noindent\hspace{2\fboxsep}\parbox{\thmw}{
  \begin{theorem}\label{#1}{\rm #2}\end{theorem}\vspace{-\lastskip}}
  \hspace{\fboxsep}}\kern-\fboxrule \vrule }}}}\newpage \hbox{%
  \lower\fboxrule\hbox{\vbox{\hbox{\vrule \kern-\fboxrule \vbox{%
  \noindent\hspace{2\fboxsep}\parbox{\thmw}{\rm #3}\hspace{\fboxsep}
  \vspace{4\fboxsep}}\kern-\fboxrule \vrule }\hrule }}}\vspace{\baselineskip}
}
\newcounter{algo}[section]
\renewcommand{\thealgo}{\thesection.\arabic{algo}}
\newcommand{\algo}[3]{\refstepcounter{algo}
\begin{center}\begin{figure}[htbp]
\framebox[\textwidth]{
\parbox{0.95\textwidth} {\vspace{\topsep}
{\bf Algorithm \thealgo : #2}\label{#1}\\
\vspace*{-\topsep} \mbox{ }\\
{#3} \vspace{\topsep} }}
\end{figure}\end{center}}
\newcommand{\calA}{{\cal A}}
\newcommand{\calD}{{\cal D}}
\newcommand{\calE}{{\cal E}}
\newcommand{\calF}{{\cal F}}
\newcommand{\calG}{{\cal G}}
\newcommand{\calH}{{\cal H}}
\newcommand{\calS}{{\cal S}}
\newcommand{\barT}{\overline{T}}
\newcommand{\barf}{\overline{f}}
\newcommand{\bart}{\overline{t}}
\newcommand{\barphi}{\overline{\phi}}
\newcommand{\barDT}{\overline{\Delta T}}
\renewcommand{\Re}{\hbox{I\hskip -2pt R}}

\newcommand{\Na}{\hbox{I\hskip -1.8pt N}}
\newcommand{\sfrac}[2]{{\scriptstyle \frac{#1}{#2}}}
\newcommand{\half}{\sfrac{1}{2}}
\newcommand{\third}{\sfrac{1}{3}}
\newcommand{\eqdef}{\stackrel{\rm def}{=}}
\newcommand{\bigfrac}[2]{\frac{\displaystyle #1}{\displaystyle #2}}
\newcommand{\bigsum}{\displaystyle \sum}
\DeclareMathOperator*\globmin{globmin}
\DeclareMathOperator*\globmax{globmax}
\newcommand{\bigglobmin}{\displaystyle \globmin}
\newcommand{\bigglobmax}{\displaystyle \globmax}
\newcommand{\ii}[1]{\{1, \ldots, #1 \}}
\newcommand{\iibe}[2]{\{ #1, \ldots, #2 \}}
\newcommand{\mystack}[2]{_{\stackrel{\scriptstyle #1}{\scriptstyle #2}}}
\newcommand{\flow}{f_{\rm low}}

\begin{document}
\title{Adaptive Regularization Algorithms with Inexact Evaluations for Nonconvex Optimization}

\author{Stefania Bellavia\footnotemark[1],
        Gianmarco Gurioli\footnotemark[2],
        Benedetta Morini\footnotemark[3], and
        Philippe L. Toint\footnotemark[4]}

\renewcommand{\thefootnote}{\fnsymbol{footnote}}
\footnotetext[1]{Dipartimento di Ingegneria Industriale,
  Universit\`{a} degli Studi, Firenze, Italy.
  Member of the INdAM Research Group GNCS.
  Email: stefania.bellavia@unifi.it.}

\renewcommand{\thefootnote}{\fnsymbol{footnote}}
\footnotetext[2]{Dipartimento di Matematica e Informatica ``Ulisse Dini",
  Universit\`{a} degli Studi, Firenze, Italy.
  Member of the INdAM Research Group GNCS.
  Email:  gianmarco.gurioli@unifi.it.}

\renewcommand{\thefootnote}{\fnsymbol{footnote}}
\footnotetext[3]{Dipartimento di Ingegneria Industriale,
  Universit\`{a} degli Studi, Firenze, Italy.
  Member of the INdAM Research Group GNCS.
  Email: benedetta.morini@unifi.it.}

\renewcommand{\thefootnote}{\fnsymbol{footnote}}
\footnotetext[4]{Namur Center for Complex Systems (naXys),
  University of Namur, 61, rue de Bruxelles, B-5000 Namur, Belgium.
  Email: philippe.toint@unamur.be.}

\date{\documentdate}

%  The manuscript

\maketitle{}

\begin{abstract}
A regularization algorithm using inexact function values and inexact
derivatives is proposed and its evaluation complexity analyzed.  This
algorithm is applicable to unconstrained problems and to problems with
inexpensive constraints (that is constraints whose evaluation and enforcement
has negligible cost) under the assumption that the derivative of highest
degree is $\beta$-H\"{o}lder continuous. It features a very flexible adaptive
mechanism for determining the inexactness which is allowed, at each iteration,
when computing objective function values and derivatives. The complexity
analysis covers arbitrary optimality order and arbitrary degree of available
approximate derivatives. It extends results of Cartis, Gould and Toint
{ [\emph{Sharp worst-case evaluation complexity bounds for arbitrary-order
nonconvex optimization with inexpensive constraints}, arXiv:1811.01220, 2018] }
on the evaluation complexity to the inexact case: if a $q$-th order minimizer
is sought using approximations to the first $p$ derivatives, it is proved that
a suitable approximate minimizer within $\epsilon$ is computed by the proposed
algorithm in at most $O\big(\epsilon^{-\frac{p+\beta}{p-q+\beta}}\big)$ iterations and
at most $O\big(|\log(\epsilon)|\epsilon^{-\frac{p+\beta}{p-q+\beta}}\big)$ approximate
evaluations.  An algorithmic variant, although more rigid in practice, can be
proved to find such an approximate minimizer in
$O\big(|\log(\epsilon)|+\epsilon^{-\frac{p+\beta}{p-q+\beta}}\big)$ evaluations.
While the proposed framework remains so far conceptual for high
degrees and orders, it is shown to yield simple and computationally realistic
inexact methods when specialized to the unconstrained and bound-constrained
first- and second-order cases. The deterministic complexity results are
finally extended to the stochastic context, yielding adaptive sample-size
rules for subsampling methods typical of machine learning.
\end{abstract}


\numsection{Introduction}

Evaluation complexity of algorithms for nonlinear and possibly nonconvex
optimization problems has been the subject of active research in recent
years. This field is concerned by deriving formal bounds on the number of
evaluations of the objective function (and possibly of its derivatives)
necessary to obtain approximate optimal solutions within a user-specified
accuracy. Until recently, the results had focused on methods using first- and
second-order derivatives of the objective function, and on convergence
guarantees to first- or second-order stationary points
\cite{Vava93,Nest04,NestPoly06,GratSartToin08,CartGoulToin11d}. Among these
contributions, \cite{NestPoly06,CartGoulToin11d} analyzed the ``regularization
method'', in which a model of the objective function around a given iterate is
constructed by adding a regularization term to the local Taylor expansion,
model which is then approximately minimized in an attempt to find a new point
with a significantly lower objective function value \cite{Grie81}. Such
methods have been shown to possess optimal evaluation complexity
\cite{CartGoulToin18a} for first- and second-order models and minimizers, and
have generated considerable interest in the research community. A
theoretically significant step was made in \cite{BirgGardMartSantToin17} for
unconstrained problems, where evaluation complexity bounds were obtained for
convergence to first-order stationary points of a simplified regularization
method using models of arbitrary degree. Even more recently,
\cite{CartGoulToin18b} proposed a conceptual unified framework subsuming all
known results for regularization methods, establishing an upper evaluation
complexity bound for arbitrary model degree and also, for the first time, for
arbitrary orders of optimality. This paper additionally covers unconstrained
problems and problems involving ``inexpensive'' constraints, that is
constraints whose evaluation/enforcement cost is negligible compared to that
of evaluating the objective function and its derivatives. It also allows for a
full range of smoothness assumptions on the objective function. Finally it
proves that the complexity results obtained are optimal in the sense that
upper and lower evaluation complexity bounds match in order. In
\cite{CartGoulToin18b}, all the above mentioned results are established for
versions of the regularization algorithms where it is assumed that objective
function values and values of its derivatives (when necessary) can be computed
exactly.

In practice, it may sometimes be difficult or impossible to obtain accurate
values of the problem's function and/or derivatives. This difficulty has been
known for a long time and has generated its own stream of results, among which
we note the trust-region method using dynamic accuracy on the objective
function and (possibly on) its gradient (see Sections~8.4.1.1 and 10.6 of
\cite{ConnGoulToin00} and \cite{BellGratRicc18}), and the purely probabilistic
approaches of \cite{PaquSche18} and \cite{BlanCartMeniSche16}. Since
unconstrained cubic regularization methods have become popular in the machine
learning community (see \cite{BottCurtNoce18} for a survey of optimization in
this area) due to their optimal complexity, several contributions have 
considered building those function and derivative's approximations by
``subsampling" the (very many) nonlinear terms whose sum defines the objective
functions typical of machine learning applications. Inexact Hessian
information is considered in
\cite{ChenJianLinZhan18,BellGuriMori18,XuRoosMaho17,XuRoosMaho18}, approximate
gradient and Hessian evaluations are used in
\cite{CartGoulToin12a,CartSche17,TripSterJinRegiJord17,YaoXuRoosMaho18},
function, gradient and Hessian values are sampled in
\cite{LiuLiuHsieTao18,BergDiouKungRoye18}. The amount of inexactness allowed
is controlled dynamically in
\cite{CartGoulToin12a,CartSche17,LiuLiuHsieTao18,ChenJianLinZhan18,BellGuriMori18}.
\vskip 5pt 

{\bf Contributions.} The present paper proposes an extension of the unifying
framework of \cite{CartGoulToin18b} for unconstrained or
inexpensively-constrained problems that allows inexact evaluations of the
objective function and of the required derivatives, in an adaptive way
inspired by the trust-region scheme of \cite[Section~10.6]{ConnGoulToin00}.
This extension has the advantage of preserving the optimal complexity of the
standard regularization methods and, as in \cite{CartGoulToin18b}, evaluation
complexity results are provided for arbitrary model degree and arbitrary order
of optimality. In particular, the proposed framework allows all combinations
of exact/inexact objective functions and derivatives of any order (including
of course degrees and orders one and two, for which simple specializations are
outlined). We also consider an interesting but practically more restrictive
variant of our algorithm for which an improved complexity can be derived. We
finally consider a stochastic version of our framework and derive rules for
sample size in the context of subsampling methods for machine learning.

The paper is organized as follows. Section~\ref{algo-s} recalls the notions of
high-order optimality proposed in \cite{CartGoulToin18b} and introduces the
general Adaptive Regularization algorithm with model of order $p$ allowing
Dynamic Accuracy (AR$p$DA). The details of how to obtain the desired relative
accuracy levels from known absolute errors are examined in
Section~\ref{rel-err-s}. The evaluation complexity of obtaining approximate
minimizers using this algorithm is then analyzed in
Section~\ref{complexity-det-s}. The algorithmic variant of the algorithm is
discussed in Section~\ref{variant-s}. The general framework is specialized to
first- and second-order optimization in Section~\ref{first-order-s}, showing
that practical implementation for low order is simple and computationally
realistic. The stochastic evaluation complexity and sampling rules for machine
learning applications are finally derived in
Section~\ref{complexity-subs-s}. Conclusions and perspectives are presented in
Section~\ref{conclusion-s}. 
\vskip 5pt 
{\bf Notations.} Unless otherwise specified, $\|\cdot\|$ denotes the standard
Euclidean norm for vectors and matrices.  For a general symmetric tensor $S$
of order $p$, we define
\beqn{Tnorm}
\|S\|_{[p]} \eqdef \max_{\|v\|=1}  | S [v]^p |
= \max_{\|v_1\|= \cdots= \|v_p\|=1} | S[v_1, \ldots, v_p] |
\eeqn
the induced Euclidean norm. We also denote by $\nabla_x^j f(x)$ the $j$-th
order derivative tensor of $f$ evaluated at $x$ and note that such a tensor is
always symmetric for any $j\geq 2$. $\nabla_x^0 f(x)$ is a synonym for $f(x)$.
$\lceil \alpha \rceil$ and $\lfloor \alpha \rfloor$ denote the
smallest integer not smaller than $\alpha$ and the largest integer not
exceeding $\alpha$, respectively. If $i$ is a non-negative integer and $\beta$ a
real in $(0,1]$ we define $(i+\beta)! = \prod_{\ell=1}^i(\ell+\beta)$. For
symmetric matrices, $\lambda_{\min}[M]$ is the leftmost eigenvalue of
$M$. $Pr[{\rm event}]$ finally denotes the probability of an event. Finally
$\globmin_{x \in \calS}f(x)$ denotes the smallest value of $f(x)$ over $x \in
\calS$.

\numsection{High-order necessary conditions and the AR$p$DA algorithm}
\label{algo-s}

Given $p \geq 1$, we consider the set-constrained optimization problem
\beqn{problem}
\min_{x \in \calF}  f(x),
\eeqn
where $\calF \subseteq \Re^n$ is closed and nonempty, and where we assume that
the \emph{values of the objective function $f$ and its derivatives must be computed
inexactly}. We also assume that $f\in  \mathcal{C}^{p,\beta}(\Re^n)$, meaning that:
\vskip 5pt 
\begin{itemize}
\item  $f$ is
$p$-times continuously differentiable, \vskip 5pt 

\item $f$ is bounded below by $f_{\rm low}$, and \vskip 5pt 

\item the $p$-th derivative tensor of $f$ at $x$ is globally H\"{o}lder
continuous, that is, there exist constants $L \geq 0$ and $\beta \in (0,1]$
such that, for all $x,y \in \Re^n$, 
\beqn{tensor-Hol}
\| \nabla_x^pf(x) - \nabla_x^pf(y) \|_{[p]} \leq L \| x-y \|^\beta.
\eeqn
\end{itemize}
The more standard case where $f$ is assumed to have Lipschitz-continuous $p$-th
derivative is recovered by setting $\beta = 1$ in the above assumptions (for
example, the choices $p=2$ and $\beta=1$ correspond to the assumption that
$f$ has a Lipschitz continuous Hessian). In what follows, we assume that
$\beta$ is known.

If we denote the $p$th degree Taylor expansion of $f$ around $x$ evaluated at $s$ by
\beqn{Taylor-def}
T_p^f(x,s) \eqdef f(x) + \sum_{\ell=1}^p \frac{1}{\ell!}\nabla_x^\ell f(x)[s]^\ell,
\eeqn
we may then define the \emph{Taylor increment} by
\beqn{deltaT-def}
\Delta T_p^f(x,s) = T_p^f(x,0) - T_p^f(x,s).
\eeqn
Under the above assumptions, we recall the crucial bounds on differences between
$f$ and its derivatives and their Taylor's expansion.

\llem{taylor-bounds-lemma}{\cite[Lemma~2.1]{CartGoulToin18b}
Let $f \in C^{p,\beta}(\Re^n)$, and $T_p^f(x,s)$ be the
Taylor approximation of $f(x+s)$ about $x$ given by \req{Taylor-def}.
Then for all $x,s \in \Re^n$,
\beqn{resf}
|f(x+s) - T_p^f(x,s)|\le  \bigfrac{L}{(p+\beta)!} \, \|s\|^{p+\beta},
\eeqn
\beqn{resder}
\| \nabla^j_x f(x+s) -  \nabla^j_s T_p^f(x,s) \|_{[j]}
\leq\bigfrac{L}{(p-j+\beta)!} \|s\|^{p+\beta-j}.
\ms (j = 1,\ldots, p).
\eeqn
}

We also follow \cite{CartGoulToin18b} and define a
$q$-th-order-necessary minimizer as a point $x \in \Re^n$
such that, for some $\delta \in (0,1]$,
\beqn{phi-def}
\phi_{f,q}^\delta(x)
\eqdef f(x)-\globmin_{\stackrel{x+d\in \calF}{\|d\|\leq\delta}}T_q^f(x,d)
= 0.
\eeqn
Observe that, in the unconstrained case, this definition subsumes the usual
optimality criteria for orders one and two, since, if $q=1$, \req{phi-def}
gives that, for any $\delta \in (0,1]$ (and in particular for $\delta = 1$),
\beqn{phi1-def}
\phi_{f,q}^\delta(x)=\|\nabla^1_xf(x)\|\delta,
\eeqn
and first-order optimality is thus equivalent to
\[
\|\nabla^1_xf(x)\|=0.
\]
Similarly, for $q=2$, \req{phi-def} is equivalent to
\beqn{derwise-exact}
\|\nabla^1_xf(x)\| = 0
\tim{ and }
\lambda_{\min}[ \nabla_x^2f(x) ] \geq 0.
\eeqn
Its properties are further discussed in \cite{CartGoulToin18b}, but we
emphasize that, for any $q\geq 1$ and in contrast with other known measures,
it varies continuously when $x$ varies continuosly in $\calF$. In the
unconstrained case, solving the global optimization problem involved in its
definition is easy for $q=1$ as the global minimizer is analytically given by
$d_* = -\delta\,\nabla_x^1f(x)/\|\nabla_x^1f(x)\|$, and also for $q=2$ using a
trust-region scheme (whose cost is essentially comparable to that of computing
the leftmost eigenvalue in \req{derwise-exact}). However this task may become
NP-hard for larger $q$. This makes $\phi_{f,q}^\delta(x)$ an essentially
theoretical tool for these cases. In any case, the computation of
$\phi_{f,q}^\delta(x)$ does not involve evaluating $f$ or any of its
derivatives, and its cost therefore does not affect the evaluation complexity
of interest here. 

If we now relax the notion of exact minimizers, we may define an
$(\epsilon,\delta)$-approximate $q$-th-order-necessary minimizer as a point
$x\in \Re^n$ 
\beqn{term-q}
\phi_{f,q}^\delta(x) \leq \epsilon \chi_q(\delta),
\eeqn
where
\beqn{chidef}
\chi_q(\delta) \eqdef \sum_{\ell=1}^q \frac{\delta^\ell}{\ell!}
\eeqn
provides a natural scaling. Again this notion reduces to familiar concepts in
the low-order unconstrained cases.  For instance, we verify that for
unconstrained problems with $q=2$, \req{term-q} requires that, if $d$ is the
global minimizer in \req{phi-def} (the solution of a trust-region problem),
\[
\max\left[0,-\left(\nabla_x^1 f(x)^Td + \half d^T\nabla_x^2f(x)d \right) \right]
\leq \epsilon (\delta+\half \delta^2),
\]
which automatically holds for any $\delta \in (0,1]$ if
$\|\nabla_x^1f(x)\| \leq \epsilon$ and 
  $\lambda_{\min}[ \nabla_x^2f(x) ] \geq -\epsilon$.
We note that, when assessing whether $x$ is an $(\epsilon,\delta)$-approximate
$q$-th-order-necessary minimizer, the global minimization in \req{phi-def} can
be stopped as soon as $\Delta T_q^f(x,d)$ exceeds $\epsilon\chi_q(\delta)$,
thereby significantly reducing the cost of this assessment.

Having defined what we mean by high-order approximate minimizers, we now turn
to des\-cribing what we mean by inaccurate objective function and derivatives
values. It is important to observe at this point that, in an optimization
problem, the role of the objective function is more central than that of any
of its derivatives, since it is the quantity we ultimately wish to decrease.
For this reason, we will handle the allowed inexactness in $f$ differently
from that in $\nabla_x^jf$: we will require an (adaptive) \emph{absolute}
accuracy for the first and a \emph{relative} accuracy for the second. In fact,
we can, in a first approach, abstract the relative accuracy requirements for
the derivatives $\nabla_x^jf(x)$ into a requirement on the relative accuracy
of $\Delta T_p^f(x,s)$.  Let $\omega \in [0,1]$ represent a relative accuracy
level and denote inexact quantities with an overbar. For what follows, we will
thus require that, if 
\beqn{batDT-def}
\barDT_p^f(x,s,\omega) = \barT_p^f(x_k,0,\omega)-\barT_p^f(x_k,s,\omega),
\eeqn
then
\beqn{barDT-acc}
|\barDT_p^f(x,s,\omega)-\Delta T_p^f(x,s)| \leq \omega \barDT_p^f(x,s,\omega).
\eeqn
It may not be obvious at this point how to enforce this relative error
bound: this is the object of Section~\ref{rel-err-s} below. For now, we
simply assume that it can be done in a finite number of evaluations of
$\{\overline{\nabla_x^jf}(x)\}_{j=1}^p$ which are inexact approximations of
$\{\nabla_x^jf(x)\}_{j=1}^p$.

Given an inexactly computed $\barDT_p^f(x,s,\omega)$ satisfying
\req{barDT-acc}, we then have to consider to compute our optimality measure
inexactly too.  Observing that the definition \req{phi-def} is independent of
$f(x)$ because of cancellation, we see that
\beqn{barphi-def}
\barphi_{f,q}^\delta(x,\omega)
= \max\Bigg[0,\globmax_{\stackrel{x+d\in \calF}{\|d\|\leq\delta}}\barDT_q^f(x,d,\omega)\Bigg].
\eeqn

Under the above assumptions, we now describe an algorithm allowing
inexact computation of both the objective function and its derivatives whose
purpose is to find (for given $q$ and a suitable relative accuracy $\omega$) a
point $x_k$ satisfying
\beqn{bar-term-q}
\barphi_{f,q}^\delta(x,\omega) \leq \frac{\epsilon}{1+\omega} \chi_q(\delta)
\eeqn
for some optimality radius $\delta \in (0,1]$.
This algorithm uses a regularized Taylor's model defined at
iteration $k$ by
\beqn{model}
m_k(s) \eqdef \barT_p^f(x_k,s,\omega_k) + \frac{\sigma_k}{(p+\beta)!} \|s\|^{p+\beta}.
\eeqn
This model is then approximately minimized and the resulting trial point is
then accepted or rejected depending on whether or not it produces a significant
decrease. This is detailed in Algorithm~\ref{algo} \vpageref{algo}.

\algo{algo}{Adaptive Regularization of order $p$ with Dynamic Accuracy (AR$p$DA)}
{\vspace*{-0.3 cm}
\begin{description}
\item[Step 0: Initialization.]
  An initial point $x_0\in\calF$  and an initial regularization parameter $\sigma_0>0$
  are given, as well as an accuracy level  $\epsilon \in (0,1)$ and an initial
  relative accuracy $\omega_0\geq 0$.  The constants  $\kappa_\omega$,
  $\delta_{-1}$, $\theta$, $\mu$, $\eta_1$, $\eta_2$, $\gamma_1$, $\gamma_2$, 
  $\gamma_3$ and $\sigma_{\min}$ are also given and satisfy
  $\theta > 0$, $\mu \in (0,1]$, $\delta_{-1} \in (0,1]$, 
  $\sigma_{\min} \in (0, \sigma_0]$,
  \beqn{eta-gamma2}
  0 < \eta_1 \leq \eta_2 < 1, \;\;
  0< \gamma_1 < 1 < \gamma_2 < \gamma_3,
  \eeqn
  \vspace*{-7mm}
  \beqn{acc-init}
  \alpha \in (0,1),
  \ms
  \kappa_\omega \in (0,\half \alpha \eta_1]
  \tim{ and }
  \omega_0 =\min\Big[\kappa_\omega,\frac{1}{\sigma_0}\Big].
  \vspace*{-3mm}
   \eeqn
  Set $k=0$.

\item[Step 1: Compute the optimality measure and check for termination. ]
  Compute $\barphi_{f,q}^{\delta_{k-1}}(x_k,\omega_k)$. If \req{bar-term-q} holds
  with $\delta = \delta_{{k-1}}$, terminate with the approximate
  solution $x_\epsilon=x_k$.
  
\item[Step 2: Step calculation. ]
 Attempt to compute a step $s_k\neq 0$ such that $x_k+s_k\in \calF$ and an
 optimality radius $\delta_k\in (0,1]$ by approximately minimizing the model
 $m_k(s)$ in the sense that 
 \vspace*{-1mm}
 \beqn{descent}
 m_k(s_k) < m_k(0)
 \vspace*{-2mm}
 \eeqn
 and
 \vspace*{-2mm}
 \beqn{mterm}
 \|s_k\| \geq \mu \epsilon^{\frac{1}{p-q+\beta}}
 \tim{ or }
 \barphi_{m_k,q}^{\delta_k}(s_k,\omega_k)
 \leq \frac{\theta\|s_k\|^{p-q+\beta}}{(p-q+\beta)!}\,\chi_q(\delta_k).
 \eeqn
 If no such step exists, terminate with the approximate solution $x_\epsilon=x_k$.
\item[Step 3: Acceptance of the trial point. ] % \ \\
 Compute $\barf_k(x_k+s_k,\omega_k)$ ensuring that
 \vspace*{-2mm}
 \beqn{Df+-DT}
 |\barf_k(x_k+s_k,\omega_k)-f(x_k+s_k)| \leq \omega_k |\barDT_p^f(x_k,s_k,\omega_k)|.
 \vspace*{-1mm}
 \eeqn
 Also ensure (by setting
 $\barf_k(x_k,\omega_k)=\barf_{k-1}(x_k,\omega_{k-1})$
 or by (re)computing $\barf_k(x_k,\omega_k)$)   that 
 \vspace*{-3mm}
 \beqn{Df-DT}
 |\barf_k(x_k,\omega_k)-f(x_k)| \leq \omega_k |\barDT_p^f(x_k,s_k,\omega_k)|.
 \vspace*{-1mm}
 \eeqn
 Then define
 \beqn{rhokdef2}
 \rho_k = \frac{\barf_k(x_k,\omega_k) - \barf_k(x_k+s_k,\omega_k)}
               {\barDT_p^f(x_k,s_k,\omega_k)}.
 \eeqn
 If $\rho_k \geq \eta_1$, then define
 $x_{k+1} = x_k + s_k$; otherwise define $x_{k+1} = x_k$.
\item[Step 4: Regularization parameter update. ]
 Set
 \vspace*{-2mm}
 \beqn{sigupdate}
 \sigma_{k+1} \in \left\{ \begin{array}{ll}
 {}[\max(\sigma_{\min}, \gamma_1\sigma_k), \sigma_k ]  & \tim{if} \rho_k \geq \eta_2, \\
 {}[\sigma_k, \gamma_2 \sigma_k ]          &\tim{if} \rho_k \in [\eta_1,\eta_2),\\
 {}[\gamma_2 \sigma_k, \gamma_3 \sigma_k ] & \tim{if} \rho_k < \eta_1.
 \end{array} \right.
 \eeqn
\item[Step 5: Relative accuracy update. ]
 Set
 \vspace*{-2mm}
 \beqn{new-acc}
 \omega_{k+1} = \min \left[ \kappa_\omega,\frac{1}{\sigma_{k+1}}\right].
 \vspace*{-2mm}
 \eeqn
 Increment $k$ by one and go to Step~1.
\end{description}
}

Some comments on this algorithm are useful at this stage.
\vskip 5pt 
\begin{enumerate}
\item
    That Step~2 may not be able, for $q>2$, to compute a nonzero step (and
    should then cause termination) can be seen by considering the following
    one-dimensional example.  Let $p=q=3$, $\calF=\Re$, $\omega_k=0$ and
    $\delta_{k-1} = 1$ and suppose that $T_3(x_k,s)= s^2-2s^3$ and also that
    $\sigma_k= 24$.  This implies that $m_k(s) = s^2-2s^3+s^4 = s^2(1-s)^2$
    and we immediately see that the origin is a global minimizer of
    $m_k(s)$. But a simple calculation shows that
    $\phi^{\delta_{k-1}}_{f,q}=T_3(x_k,0)-T_3(x_k,1) = 1$ and hence
    termination will not occur in Step~1 if $\epsilon <1/\chi_3(1) = 4/7$.
    As a consequence, as was pointed out in
    \cite{CartGoulToin18b}, the possibility of a zero $s_k$ cannot be ignored
    in Step~2. In this case, it is not possible to satisfy \eqref{descent}
    and the algorithm terminates with $x_{\varepsilon}=x_k$. 
    It has been proved in  \cite[Lemma 2.6]{CartGoulToin18b} that this is
    acceptable (see also Lemma \ref{zerostep} below).
    \vskip 5pt 
\item Our assumption \req{barDT-acc} is used three times in the algorithm: in
    Step~1 for computing $\barphi_{f,q}^{\delta_{k-1}}(x_k,\omega_k)$ and in
    Step~2 when computing $s_k$ and $\barphi_{m_k,q}^{\delta_k}(s_k,\omega_k)$.\vskip 5pt 
\item As indicated above, we require a bound on the absolute error in the
    objective function value: this is the object of \req{Df+-DT} and
    \req{Df-DT}, where we introduced the notation $\barf_k(x_k,\omega_k)$ to
    denote an inexact approximation of $f(x_k)$. Note that a new value of
    $\barf_k(x_k,\omega_k)$ should be computed to ensure \req{Df-DT} in Step~3
    only if $k > 0$ and
    $\omega_{k-1}\barDT_p^f(x_{k-1},s_{k-1},\omega_{k-1})>
      \omega_k\barDT_p^f(x_k,s_k,\omega_k)$.   
    If this is the case the (inexact) function value is computed twice per
    iteration instead of just once.\vskip 5pt 
\item At variance with the trust-region method with dynamic accuracy of
    \cite[Section~10.6]{ConnGoulToin00} and \cite{BellGratRicc18}, we do not
    recompute approximate values of the objective function at $x_k$ once the
    computation of $s_k$ is complete (provided we can ensure \req{barDT-acc},
    as discussed in Section~\ref{rel-err-s}).\vskip 5pt 
\item If $\|s_k\|\geq \mu \epsilon^{\frac{1}{p-q+\beta}}$ in Step~2, then the
    (potentially costly) calculation of
    $\barphi_{m_k,q}^{\delta_k}(s_k,\omega_k)$ is unecessary and $\delta_k$
    may be chosen arbitrarily in $(0,1]$.\vskip 5pt 
\item We call iteration $k$ \emph{successful} when $\rho_k \geq \eta_1$ and
    $x_{k+1}= x_k+s_k$.  The iteration is called \emph{unsuccessful}
    otherwise, and $x_{k+1}= x_k$ in this case. We use the notation
    \beqn{calSk-def}
    \calS_k = \{ j \in \iibe{0}{k} \mid \rho_j \geq \eta_1 \}
    \eeqn
    to denote the set of successful iterations of index at most $k$.\vskip 5pt 
\item As indicated above, ensuring \req{barDT-acc} may require a certain
    number of (approximate) evaluations of the derivatives of $f$.  For a
    single iteration of the algorithm, these evaluations are always at the
    current iterate $x_k$.
    \item It is worth noting that from \req{eta-gamma2}, \req{acc-init},
      \req{sigupdate} and \req{new-acc},   together with the positivity
      of $\sigma_0$ and $\sigma_{\min}$,
    \beqn{omega_bound}
    0<\omega_k\le \kappa_\omega < 1.
    \eeqn
\end{enumerate}
\clearpage

We now state some properties of the algorithm that are derived without
modification from the case where the computation of $f$ and its derivatives
are exact.

\llem{SvsU}{\cite[Theorem~2.1]{CartGoulToin11d}
The mechanism of the AR$p$DA algorithm ensures that, if
\vspace*{-2mm}
\beqn{sigmax}
\sigma_{k} \leq \sigma_{\max},
\vspace*{-2mm}
\eeqn
for some $\sigma_{\max} > 0$, then
\vspace*{-2mm}
\beqn{unsucc-neg}
k +1 \leq |\calS_k| \left(1+\frac{|\log\gamma_1|}{\log\gamma_2}\right)+
\frac{1}{\log\gamma_2}\log\left(\frac{\sigma_{\max}}{\sigma_0}\right).
\eeqn
}

\noindent
This shows that the number of unsuccessful iterations must remain a fixed
proportion of that of the successful ones.

\llem{step-ok-l}{\cite[Lemma~2.5]{CartGoulToin18b}
Suppose that $s_k^*\neq 0$ is a global minimizer of $m_k(s)$
under the constraint that $x_k+s \in \calF$, such $m_k(s_k^*)< m_k(0)$.
Then there exist a neighbourhood of $s_k^*$ and a range of sufficiently
small $\delta$ such that \req{descent} and the second part of \req{mterm}
hold for any $s_k$ in the intersection of this neighbourhood with $\calF$
and any $\delta_k$ in this range.
}

\noindent
This last lemma thus ensures that the algorithm is well-defined when
$s_k\neq0$. The lemma below shows that it is reasonable to
terminate the algorithm whenever  a nonzero descent step cannot be
computed.

\llem{zerostep}{\cite[Lemma~2.6]{CartGoulToin18b}
Suppose that the algorithm terminates in Step~2 of iteration $k$ with 
$x_{\varepsilon}=x_k$.
Then there exists a $\delta \in  (0, 1]$  such that \req{bar-term-q} holds  for $x = x_{\varepsilon}$.
}

\numsection{Enforcing the relative error on Taylor increments}
\label{rel-err-s}

We now return to the question of enforcing \req{barDT-acc}. For improved
readability, we temporarily ignore the iteration index $k$.

\subsection{The accuracy checks}

While there may be circumstances where \req{barDT-acc} can be enforced
directly, we consider here that the only control the user has on the accuracy
of $\barDT_p^f(x,s,\omega)$ is by enforcing bounds $\{\varepsilon_j\}_{j=1}^p$
on the absolute errors on the derivative tensors $\{\nabla_x^jf(x)\}_{j=1}^p$.
In other words, we seek to ensure \req{barDT-acc} by selecting
absolute accuracies $\{\varepsilon_j\}_{j=1}^p$ such that, when
\beqn{vareps-j}
\|\overline{\nabla_x^jf}(x)-\nabla_x^jf(x)\|_{[j]} \leq \varepsilon_j
\tim{for} j \in \ii{p},
\eeqn
the desired accuracy requirement follows. 

In all cases described below, the process can be viewed as an iteration with
four main steps.  The first is to compute the relevant approximate derivative
satisfying \req{vareps-j} for given values of $\{\varepsilon_j\}_{j=1}^p$.
The second is to use these approximate derivatives to compute the desired
Taylor increment and associated quantities.  Tests are then performed in the
third step to verify the desired accuracy requirements and terminate if they
are met. If not the case, the absolute accuracies $\{\varepsilon_j\}_{j=1}^p$
are then decreased before a new iteration is started.

As can be expected, a suitable relative accuracy requirement will be
achievable as long as $\barDT_p^f(x,s,\omega)$ remains safely away from zero,
but, if exact computations are to be avoided, we may have to accept a simpler
absolute accuracy guarantee when $\barDT_p^f(x,s,\omega)$ vanishes.

We then formalize the resulting accuracy tests in the VERIFY algorithm, stated as
Algorithm~\ref{verify} \vpageref{verify}. 

Assume that for a vector $v_\omega$, a bound $\delta \geq \|v_\omega\|$, a
degree $r$, the requested relative and absolute accuracies $\omega$ and
$\xi>0$, the increment $\barDT_r(x,v_\omega,\omega)$ are given.  We intend to
use the algorithm for $\barDT_q^f(x,v_\omega,\omega)$,
$\barDT_p^f(x,v_\omega,\omega)$ and $\barDT_q^{m_k}(x,v_\omega,\omega)$.  For
keeping our development general, we use the notations
$\barDT_r(x,v_\omega,\omega)$ and $\Delta T_r(x,v_\omega)$ without
superscript.  Moreover, we assume that the current absolute accuracies
$\{\zeta_j\}_{j=1}^r$ of the derivatives of $\overline{
  T}_r(x,v_\omega,\omega)$ with respect to $v_\omega$ at $v_\omega=0$ are
given. Because it will be the case below, we assume for simplicity that
$\barDT_r(x,v_\omega,\omega) \geq 0$.

\vspace*{-3mm}
\algo{verify}{Verify the accuracy of $\barDT_r(x,v_\omega,\omega)$\protect\\
\centerline{
${\tt flag} = {\rm VERIFY}
\Big(\delta,\barDT_r(x,v_\omega,\omega),\{\zeta_j\}_{j=1}^{r},\omega,\xi\Big)$
}}
{
\parbox[t]{13.5cm}{
  \vspace*{-5mm}
Set ${\tt flag} = 0$.
\begin{itemize}
\item If
  \beqn{verif-term-1}
  \barDT_r(x,v_\omega,\omega) = 0
  \tim{ and }
  \max_{j\in \ii{r}}\zeta_j \leq \xi,
  \eeqn
  set {\tt flag} = 1.
\item Else, if
  \vspace*{-3mm}
  \beqn{verif-term-2}
  \barDT_r(x,v_\omega,\omega) > 0
  \tim{ and }
  \sum_{j=1}^r \frac{\zeta_j}{j!}\delta^j \leq \omega \barDT_r(x,v_\omega,\omega),
  \vspace*{-3mm}
  \eeqn
  set {\tt flag} = 2.
\item Else, if
  \vspace*{-3mm}
  \beqn{verif-term-3}
  \barDT_r(x,v_\omega,\omega) > 0
  \tim{ and }
  \sum_{j=1}^r \frac{\zeta_j}{j!}\delta^j \leq \xi \chi_r(\delta),
  \vspace*{-3mm}
  \eeqn
  set {\tt flag} = 3.
  \end{itemize}
  }
}

\noindent
Let us now consider what properties are ensured for the various possible
values of ${\tt flag}$.

\llem{verify-l}{
Suppose that
\beqn{rvareps-j}
\Big\|\Big[\overline{\nabla_{v_\omega} ^j T_r}(x,v_\omega)\Big]_{v_\omega=0}
-\Big [\nabla_{v_\omega} ^j T_r(x,v_\omega)\Big]_{v_\omega=0}\Big\|_{[j]}
\leq \zeta_j \tim{for} j \in \ii{r}
\eeqn
and $\omega\in(0,1)$. Then we have that
\begin{itemize}
\item if
  \vspace*{-3mm}
  \beqn{max-abs-acc}
    \max_{j \in \ii{r}}\zeta_j \leq \xi,
  \eeqn
  then the VERIFY algorithm returns a nonzero ${\tt flag}$,
\item if the VERIFY algorithm terminates with ${\tt flag} = 1$, then
  $\barDT_r(x,v_\omega,\omega)= 0$ and
  \beqn{verif-prop-1}
  \left|\barDT_r(x,v,\omega)- \Delta T_r(x,v)\right|\leq \xi \chi_r(\|v\|)
  \ms\tim{for all $v$,}
  \eeqn
\item if the VERIFY algorithm terminates with ${\tt flag} = 2$, then
  $\barDT_r(x,v_\omega,\omega)>0$ and
  \beqn{verif-prop-2}
  \left|\barDT_r(x,v,\omega)- \Delta T_r(x,v)\right|
  \leq \omega \barDT_r(x,v_\omega,\omega),
  \ms\tim{for all $v$ with $\|v\|\leq \delta$,} 
  \eeqn
\item if the VERIFY algorithm terminates with ${\tt flag} = 3$, then
  $\barDT_r(x,v_\omega,\omega)>0$ and
  \beqn{verif-prop-3}
  \max\left[\barDT_r(x,v_\omega,\omega),
    \left|\barDT_r(x,v,\omega)- \Delta T_r(x,v)\right|\right]
  \leq \frac{\xi}{\omega} \chi_r(\delta)
  \tim{for all $v$ with $\|v\|\leq \delta$.}
  \eeqn
\end{itemize}
}

\begin{proof}
We first prove the first proposition. If $\barDT_r(x,v_\omega,\omega) = 0$ and
\req{max-abs-acc}, then \req{verif-term-1} ensures that ${\tt flag = 1}$ is
returned. If $\barDT_r(x,v_\omega,\omega) > 0$, from \req{chidef} and
\req{max-abs-acc}, we deduce that
\[
\sum_{j=1}^r \frac{\zeta_j}{j!}\delta^j
\leq \Big[\max_{j\in\ii{r}}\zeta_j \Big] \chi_r(\delta)
\leq \xi \chi_r(\delta)
\]
also causing termination with ${\tt flag} = 3$ because of \req{verif-term-3}
if it has not occurred with ${\tt flag} = 2$ because of \req{verif-term-2},
hence proving the first proposition.
 
Consider now the three possible termination cases and suppose first that
termination occurs with {\tt flag} = 1.  Then, using the triangle inequality,
\req{rvareps-j}, \req{verif-term-1} and \req{chidef}, we have that, 
for any $v$,
\vspace*{-2mm}
\[
\left|\barDT_r(x,v,\omega)- \Delta T_r(x,v)\right|
\leq \sum_{j=1}^r \frac{\zeta_j}{j!} \|v\|^j
\leq \xi \chi_r(\|v\|)
\vspace*{-2mm}
\]
yielding \req{verif-prop-1}. Suppose now that {\tt flag = 2}.
Then \req{verif-term-2} holds and for any $v$ with
$\|v\|\leq \delta$,
\vspace*{-2mm}
\[
\left|\barDT_r(x,v,\omega)- \Delta T_r(x,v)\right|
\leq \sum_{j=1}^r \frac{\zeta_j}{j!} \|v\|^j
\leq \sum_{j=1}^r \frac{\zeta_j}{j!} \delta^j
\leq \omega \barDT_r(x,v_\omega,\omega),
\vspace*{-2mm}
\]
which is \req{verif-prop-2}. Suppose finally that {\tt flag = 3}. Since
termination did not occur in \req{verif-term-2}, we have that
\beqn{small-barDT}
0< \omega \barDT_r(x,v_\omega,\omega)\leq \xi \chi_r(\delta).
\eeqn
Furthermore, \req{verif-term-3} implies that, for any $v$ with
$\|v\|\leq\delta$, 
\[
\left|\barDT_r(x,v,\omega)- \Delta T_r(x,v)\right|
\leq \sum_{j=1}^r \frac{\zeta_j}{j!} \|v\|^j
\leq \sum_{j=1}^r \frac{\zeta_j}{j!} \delta^j
\leq \frac{\xi}{\omega} \chi_r(\delta).
\]
This inequality and \req{small-barDT} together imply \req{verif-prop-3}.
\end{proof}%epr

\noindent
Clearly, the outcome corresponding to our initial aim to obtain a relative
error at most $\omega$ corresponds to the case where ${\tt flag} = 2$. As we
will see below, the two other cases are also useful.

\subsection{Computing $\barphi_{f,q}^{\delta_{k-1}}(x_k,\omega_k)$}

\noindent
We now consider, in Algorithm~\ref{step1}, how to compute
the optimality measure $\barphi_{f,q}^{\delta_{k-1}}(x_k,\omega_k)$ in Step~1
of the AR$p$DA algorithm.

\algo{step1}{Modified Step 1 of the AR$p$DA algorithm}{
\vspace*{-5mm}
\begin{description}
\item[Step 1: Compute the optimality measure and check for termination.]
\mbox{}\\
\begin{description}
\item[Step 1.0: ] 
     The iterate $x_k$ and the radius $\delta_{k-1}\in (0,1]$ are
     given, as well as constants $\gamma_\varepsilon \in (0,1)$ and
     $\kappa_\varepsilon > 0$. Set $i_\varepsilon=0$.
\item[Step 1.1: ] 
     Choose an initial set of derivative absolute accuracies
     $\{\varepsilon_{j,0}\}_{j=1}^p$ such that
     \beqn{upp-vareps}
     \varepsilon_{j,0} \leq \kappa_\varepsilon \tim{ for } j \in \ii{p}.
     \eeqn
\vspace*{-2mm}
\item[Step 1.2: ] If unavailable, compute
  $\{\overline{\nabla_x^jf}(x_k)\}_{j=1}^q$ satisfying
  \[
  \|\overline{\nabla_x^jf}(x)-\nabla_x^jf(x)\|_{[j]} \leq \varepsilon_{j,i_\varepsilon}
\tim{for} j \in \ii{q}.
  \]
\item[Step 1.3: ] Solve
\vspace*{-2mm}
\[
\bigglobmax_{\stackrel{x_k+d\in \calF}{\|d\|\leq\delta_{k-1}}}\barDT_q^f(x_k,d,\omega_k),
\vspace*{-1mm}
\]
to obtain the maximizer $d_k$ and the corresponding Taylor increment
$\barDT_q^f(x_k,d_k,\omega_k)$. Compute
\vspace*{-3mm}
\[
{\tt flag} =
{\rm VERIFY}\Big(\delta_{k-1},\barDT_q^f(x_k,d_k,\omega_k),\{\epsilon_j\}_{j=1}^q,
                 \omega_k,\half\omega_k\epsilon \Big).
\vspace*{-2mm}
\]
\item[Step 1.4: ] 
Terminate the AR$p$DA algorithm with the approximate solution $x_\epsilon=x_k$
if ${\tt flag} = 1$, or if ${\tt flag} = 3$, or if ${\tt flag} = 2$ and
\req{bar-term-q} holds with $\delta = \delta_{k-1}$. Also go to Step~2 of the
AR$p$DA algorithm if  ${\tt flag} = 2$ but \req{bar-term-q} fails.
\item[Step 1.5: ] Otherwise (i.e.\ if ${\tt flag} = 0$), set
  \vspace*{-2mm}
  \beqn{decrease-vareps}
  \varepsilon_{j,i_\varepsilon+1} = \gamma_\varepsilon\varepsilon_{j,i_\varepsilon}
  \tim{ for } j \in \ii{p},
  \vspace*{-2mm}
  \eeqn
  increment $i_\varepsilon$ by one and return to Step~1.1.
\end{description}
\end{description}
}

\noindent
We immediately observe that Algorithm~\ref{step1} terminates in a finite
number of iterations, since it does so as soon as ${\tt flag} > 0$, which,
because of the first proposition of Lemma~\ref{verify-l}, must happen after a
finite number of passes in iterations using \req{decrease-vareps}. We discuss
in Section~\ref{compexity-single-s} exactly how many such decreases might be
needed.

We now verify that terminating the AR$p$DA algorithm as indicated in this
modified version of Step~1 provides the required result. We start noting that,
if $x_k$ is an isolated feasible point (i.e. such that the intersection of any
ball of radius $\delta_{k-1}>0$ centered at $x_k$ with $\calF$ is reduced to
$x_k$), then clearly $d_k=0$ and thus, irrespective of $\omega_k$ and
$\delta_{k-1} > 0$,
\beqn{phi-isolated}
\phi_{f,q}^{\delta_{k-1}}(x_k) = 0 = \barDT_q^f(x_k,d_k,\omega_k)
= \barphi_{f,q}^{\delta_{k-1}}(x_k,\omega_k),
\eeqn
which means that $\barphi_{f,q}^{\delta_{k-1}}(x_k,\omega_k)$ is a faithful indicator
of optimality at $x_k$.

\llem{acc-S1-l}{
If the AR$p$DA algorithm terminates within Step~1.4, then
\beqn{term-q-k}
\phi_{f,q}^{\delta_{k-1}}(x_k) \leq \epsilon \chi_q(\delta_{k-1})
\eeqn
and $x_k$ is a $(\epsilon,\delta_{k-1})$-approximate
$q$-th-order-necessary minimizer. Otherwise Algorithm~\ref{step1} terminates with
\beqn{good-phi}
(1-\omega_k) \barphi_{f,q}^{\delta_{k-1}}(x_k,\omega_k)
\leq \phi_{f,q}^{\delta_{k-1}}(x_k)
\leq (1+\omega_k) \barphi_{f,q}^{\delta_{k-1}}(x_k,\omega_k).
\eeqn
}

\begin{proof}
We first notice that Step~1.2 of Algorithm~\ref{step1} yields
\req{rvareps-j} with $T_r=T_r^{f}$, $r=q$ and
$\{\zeta_j\}_{j=1}^{r}=\{\varepsilon_{j,i_\varepsilon}\}_{j=1}^q$. Furthermore,
$\omega=\omega_k\in(0,1)$, so that the assumptions of Lemma~\ref{verify-l} are
satisfied. If $x_k$ is an isolated feasible point, the lemma's conclusions 
directly follow from \req{phi-isolated}.  Assume therefore that $x_k$ is not
an isolated feasible point and note first that, because Step~1.3 finds the
global maximum of $\barDT_q^f(x_k,d,\omega_k)$, we have that
$\barDT_q^f(x_k,d_k,\omega_k)\geq 0$. Suppose now that, in Step~1.3, the
VERIFY algorithm returns {\tt flag} = 1 and thus that
$\barDT_q^f(x_k,d_k,\omega_k)=0$. This means that $x_k$ is a global
minimizer of $\barT_q^f(x_k,d,\omega_k)$ in the intersection of a ball of
radius $\delta_{k-1}$ and $\calF$ and $\barDT_q^f(x_k,d,\omega_k)\leq 0$ for
any $d$ in this intersection. Thus, for any such $d$, we obtain from
\req{verif-prop-1} with $\xi = \half \omega_k \epsilon$ that
\[
\Delta T_q^f(x_k,d)
\leq  \barDT_q^f(x_k,d,\omega_k) + 
      \left|\barDT_q^f(x_k,d,\omega_k)- \Delta T_q^f(x_k,d)\right| 
\leq \half \omega_k \epsilon \chi_q(\delta_{k-1}),
\]
which, since $\omega_k\leq 1$, implies \req{term-q-k}.  Suppose next that the
VERIFY algorithm returns {\tt flag} = 3. Then $\barDT_q^f(x_k,d_k,\omega_k)>0$
and thus $d_k\neq 0$. Using the fact that the nature of Step~1.3 ensures that
$\barDT_q^f(x_k,d,\omega_k) \leq \barDT_q^f(x_k,d_k,\omega_k)$ for $d$ with
$\|d\|\leq\delta_{k-1}$ we have, using \req{verif-prop-3} with $\xi = \half
\omega_k \epsilon$, that, for all such $d$,
\[
  \begin{array}{lcl}
  \Delta T_q^f(x_k,d)
  & \leq & \barDT_q^f(x_k,d,\omega_k) + 
  \left|\barDT_q^f(x_k,d,\omega_k)- \Delta T_q^f(x_k,d)\right| \\*[1ex]
  & \leq & \barDT_q^f(x_k,d_k,\omega_k) + 
  \left|\barDT_q^f(x_k,d,\omega_k)- \Delta T_q^f(x_k,d)\right| \\*[1ex]
  & \leq & \epsilon \chi_q(\delta_{k-1})
  \end{array}
\]
yielding  \req{term-q-k}. If the VERIFY algorithm returns {\tt flag = 2},
then, for any $d$ with $\|d\|\leq \delta_{k-1}$,
\[
  \Delta T_q^f(x_k,d)
  \leq \barDT_q^f(x_k,d,\omega_k) + 
  \left|\barDT_q^f(x_k,d,\omega_k)- \Delta T_q^f(x_k,d)\right| 
  \leq (1+\omega_k)\barDT_q^f(x_k,d_k,\omega_k).
\]
Thus, for all $d$ with $\|d\|\leq \delta_{k-1}$,
\beqn{good-phi-1}
  \max\Big[0,\Delta T_q^f(x_k,d)\Big]
  \leq (1+\omega_k)\max\Big[0,\barDT_q^f(x_k,d_k,\omega_k)\Big]
  = (1+\omega_k)\barphi_{f,q}^{\delta_{k-1}}(x_k,\omega_k).
\eeqn
But termination implies that \req{bar-term-q} holds for $\delta=
\delta_{k-1}$, and \req{term-q-k} follows with this value of $\delta$.
Finally, if the AR$p$DA algorithm does not terminates within
Step~1.4 but Algorithm~\ref{step1} terminates, it must be because the VERIFY
algorithm returns {\tt flag} = 2.  This implies, as above, that
\req{good-phi-1} holds, which is the rightmost part of
\req{good-phi}. Similarly, for any $d$ with $\|d\|\leq \delta_{k-1}$,
\[
  \begin{array}{lcl}
  \Delta T_q^f(x_k,d)
  & \geq & \barDT_q^f(x_k,d,\omega_k) - 
  \left|\barDT_q^f(x_k,d,\omega_k)- \Delta T_q^f(x_k,d)\right| \\*[1ex]
  & \geq & \barDT_q^f(x_k,d,\omega_k) - \omega_k\barDT_q^f(x_k,d_k,\omega_k).
  \end{array}
\]
\vspace*{-2mm}
Hence
\[
  \begin{array}{lcl}
  \bigglobmax_{\mystack{x_k+d \in \calF}{\|d\|\leq \delta_{k-1}}}\Delta T_q^f(x_k,d)
  & \geq & \bigglobmax_{\mystack{x_k+d \in \calF}{\|d\|\leq \delta_{k-1}}}
  \left[ \barDT_q^f(x_k,d,\omega_k) -
     \omega_k\barDT_q^f(x_k,d_k,\omega_k)\right] \\*[2ex]
  & = & (1-\omega_k)\barDT_q^f(x_k,d_k,\omega_k).
  \end{array}
\]
Since $\barDT_q^f(x_k,d_k,\omega_k)>0$ when the VERIFY algorithm returns
{\tt flag} = 2, we then obtain that, for all $\|d\|\leq \delta_{k-1}$,
\[
  \max\Big[0,\globmax_{\mystack{x_k+d \in \calF}{\|d\|\leq \delta_{k-1}}}\Delta T_q^f(x_k,d)\Big]
  \geq \max\Big[0,(1-\omega_k)\barDT_q^f(x_k,d_k,\omega_k)\Big]
  = (1-\omega_k)\barphi_{f,q}^{\delta_{k-1}}(x_k,\omega_k),
\]
which is the leftmost part of \req{good-phi}.
\end{proof}%epr

\subsection{Computing $s_k$}

We now consider computing $s_k$ at Step~2 of the AR$p$DA algorithm. The
process is more complicated than for Step~1, as it potentially involves two
situations in which one wishes to guarantee a suitable relative error. The
first is when minimizing the model
\[
m_k(s) = \barf(x_k,\omega_k) -\barDT_p^f(x_k,s,\omega_k)
         + \frac{\sigma_k}{(p+\beta)!}\|s\|^{p+\beta}
\]
or, equivalently, maximizing
\beqn{minus-mk}
-m_k(s) =  - \barf(x_k,\omega_k)+\barDT_p^f(x_k,s,\omega_k)
           - \frac{\sigma_k}{(p+\beta)!}\|s\|^{p+\beta},
\eeqn
and the second is when globally minimizing the model's Taylor expansion taken
at $x_k+s_k$ in a neighbourhood of diameter $\delta_k$.  The first of these
situations can be handled in a way very similar to that used above for
computing $\barphi_{f,q}^{\delta_{k-1}}(x_k)$ in Step~1: given a set of
approximate derivatives, a step $s_k$ is computed such that it satisfies
\req{descent} and \req{mterm}, the relative error of the associated
$\barDT_p^f(x_k,s_k,\omega_k)$ is then evaluated and, if it is insufficient,
the accuracy on the derivative approximations improved and the process
restarted. If the relative error on $\barDT_p^f(x_k,s_k,\omega_k)$ is
satisfactory and the first test of \req{mterm} fails, it remains to check that
the relative error on $\barphi_{m_k,q}^{\delta_k}(s_k,\omega_k)$ is also
satisfactory.  
Moreover, as in the original AR$p$DA algorithm, we have to take
into account the possibility that minimizing the model might result in a
vanishing decrease. The resulting somewhat involved process is formalized in
Algorithm~\ref{step2} \vpageref{step2}.

\algo{step2}{Modified Step 2 of the AR$p$DA algorithm}{
\vspace*{-5mm}
\begin{description}
\item[Step 2: Step calculation. ]
\mbox{}\\
%\vspace*{-6mm}
\begin{description}
\item[Step 2.0: ] The iterate $x_k$, the radius $\delta_{k-1}\in (0,1]$, the
constants $\gamma_\varepsilon \in (0,1)$, $\vartheta\in(0,1)$,
the counter $i_\varepsilon$ and the absolute accuracies
$\{\varepsilon_{j,i_\varepsilon}\}_{j=1}^p$ are given. 
\item[Step 2.1: ] If unavailable, compute
  $\{\overline{\nabla_x^jf}(x_k)\}_{j=1}^p$ satisfying \req{vareps-j}
  with $\varepsilon_j= \varepsilon_{j,i_\varepsilon}$ for $j \in \ii{p}$.
\item[Step 2.2: ]
    \begin{itemize}
    \item Attempt to compute a step $s_k\neq 0$ with $x_k+s_k\in \calF$ such that
  \req{descent} holds.
  \item If this not possible, set ${\tt flag}_s=1$ and go to Step~2.3.
  \item Otherwise, pursue the approximate
  minimization of the model $m_k(s)$ for $x_k+s_k\in \calF$ in order to
  satisfy \req{mterm}, yielding a step  $s_k$, a decrease
  $\barDT_p^f(x_k,s_k,\omega_k)$ and, if the first part of \req{mterm} fails,
  the global maximizer $d_k^{m_k}$ of
  $\barDT_q^{m_k}(s_k,d,\omega_k)$ subject to $\|d\|\leq \delta_k$
  and $x_k+s_k+d\in \calF$, together  with the corresponding Taylor increment
  $\barDT_q^{m_k}(s_k,d_k^{m_k},\omega_k)$.
  \item Compute
  \vspace*{-3mm}
  \[
  {\tt flag}_s = {\rm VERIFY}\Big(
  \|s_k\|,\barDT_p^f(x_k,s_k,\omega_k),\{\epsilon_j\}_{j=1}^p,\omega_k,\half
  \omega_k\epsilon \Big).
  \vspace*{-2mm}
  \]
  If ${\tt flag}_s = 0$ go to Step~2.5.
  \end{itemize}
\item[Step 2.3: ] 
  If ${\tt flag}_s = 1$ or ${\tt flag}_s= 3$, compute
  \vspace*{-3mm}
  \[
  \bigglobmin_{x_k+s\in \calF} m_k(s),
  \vspace*{-3mm}
  \]
  to obtain the minimizer $s_k$, $\barDT_p^f(x_k,s_k,\omega_k)$.\\
  Set $d_k^{m_k}=0=\barDT_q^{m_k}(s_k,d_k^{m_k},\omega_k)$ and compute
  \vspace*{-3mm}
  \[
  {\tt flag}_s = {\rm VERIFY}\Big(
  \|s_k\|,\barDT_p^f(x_k,s_k,\omega_k),\{\epsilon_j\}_{j=1}^p,\omega_k,\half
  \omega_k\epsilon \Big).
  \vspace*{-2mm}
  \]
  If ${\tt flag}_s = 0$ go to Step~2.5.
\item[Step 2.4: ]
  If ${\tt flag}_s = 1$ or ${\tt flag}_s= 3$, terminate the AR$p$DA algorithm with
  $x_\epsilon=x_k$.  Otherwise, if $\|s_k\| \geq \mu
  \epsilon^{\frac{1}{p-q+\beta}}$ or if $\|s_k\| < \mu
  \epsilon^{\frac{1}{p-q+\beta}}$ and
  \[
  \hspace*{-8mm}{\tt flag}_d \!=\! {\rm VERIFY}\Big(
  \delta_k,\barDT_q^{m_k}\!(s_k,d_k^{m_k}\!,\omega_k),\{3\epsilon_j\}_{j=1}^q,\omega_k,
  \frac{\vartheta(1-\kappa_\omega)}{(1+\kappa_\omega)^2}\, \frac{\omega_k\epsilon}{2} \Big) > 0,
  \vspace*{-2mm}
  \]
  go to Step~3 of the AR$p$DA algorithm with the step $s_k$, the
  associated $\barDT_p^f(x_k,s_k,\omega_k)$ and $\delta_k$. 
\item[Step 2.5: ] Set (if ${\tt flag}_s =0$ or ${\tt flag}_d =0$),
  \beqn{decrease-vareps-2}
  \varepsilon_{j,i_\epsilon+1} = \gamma_\varepsilon\varepsilon_{j,i_\epsilon}
  \tim{ for } j \in \ii{p},
  \eeqn
  increment $i_\varepsilon$ by one and go to Step~2.1.
%  }
\end{description}
\end{description}
}

\noindent
Observe that, in Step~2.2, $d_k^{m_k}$ and $\barDT_q^{m_k}(s_k,d_k^{m_k},\omega_k)$
result from the computation of
$\barphi_{m_k,q}^{\delta_k}(s_k,\omega_k)$ which is necessary to verify
the second part of \req{mterm}. Note also that we have specified,
in the call to VERIFY in Step~2.4 of Algorithm~\ref{step2}, absolute
accuracy values equal to $\{3\epsilon_j\}_{j=1}^q$. This is because
this call aims at checking  the accuracy of the Taylor expansion of the
\emph{model} and  the derivatives which are then approximated are
not $\{\overline{\nabla_x^j f}(x_k)\}_{j=1}^q$, but $\{\overline{\nabla_d^j
T}{\,}_q^{m_k}(s_k,0)\}_{j=1}^q$.  It is easy to verify that these
(approximate) derivatives are given by 
\beqn{ders-mk-plus}
\overline{\nabla_d^jT}{\,}_q^{m_k}(s_k,0)
=\sum_{\ell=j}^p \frac{\overline{\nabla_x^\ell f}(x_k)\|s_k\|^{\ell-j}}{(\ell-j)!}
+ \left[\nabla_s^j\|s\|^{p+\beta}\right]_{s=s_k},
\eeqn
where the last term of the right-hand side is exact.  This yields the
following error bound. 

\llem{abserr-derm-l}{
Suppose that $\|s_k\| \leq \mu\epsilon^{\frac{1}{p-q+\beta}}$. Then, for all $j \in \ii{p}$,
\beqn{abserr-derm}
\left|\overline{\nabla_d^jT}{\,}_q^{m_k}(s_k,0)-\nabla_d^jT_q^{m_k}(s_k,0)\right|
\leq 3 \,\varepsilon_j.
\eeqn
}

\begin{proof}
Using the triangle inequality, \req{ders-mk-plus}, the inequality
$\|s_k\| \leq \mu\epsilon^{\frac{1}{p-q+\beta}}\leq \mu$ and \req{chidef}, we
have that, for all $j\in\ii{p}$,
\[
\left|\overline{\nabla_d^j T}{\,}_q^{m_k}(s_k,0)-\nabla_d^jT_q^{m_k}(s_k,0)\right|
\leq \sum_{\ell=j}^p \frac{\varepsilon_j\|s_k\|^{\ell-j}}{(\ell-j)!}
\leq \varepsilon_j \sum_{\ell=j}^p \frac{\mu^{\ell-j}}{(\ell-j)!}
\leq \varepsilon_j (1+\chi_p(\mu))
\]
and \req{abserr-derm} follows since $\chi_p(\mu) \leq 2 \mu$.
\end{proof}\\%epr

Again, Algorithm~\ref{step2} must terminate in a finite number of iterations.
Indeed, if after finitely many iterations ${\tt flag}_s = 1$ or ${\tt flag}_s
= 3$ at the start of Step~2.4, the conclusion is obvious. Suppose now that
${\tt flag}_s = 2$ at all iterations. If $\|s_k\| <
\mu\epsilon^{\frac{1}{p-q+\beta}}$ always hold, the first proposition of
Lemma~\ref{verify-l} ensures that ${\tt flag}_d > 0$ after finitely many
decreases in \req{decrease-vareps-2}, also causing termination. Termination
might of course occur if $\|s_k\| \geq \mu\epsilon^{\frac{1}{p-q+\beta}}$
before this limit.

The next Lemma characterizes the outcomes of Algorithm~\ref{step2}.

\llem{good-step2}{Suppose that the modified Step~2
  is used in the AR$p$DA algorithm. If this algorithm terminates
  within that step, then there exists a $\delta \in  (0, 1]$  such that
  \req{bar-term-q} holds for $x = x_{\varepsilon}$ or
  \beqn{strong-opt}
  \phi_{f,p}^{\|s_k\|}(x_k) \leq \epsilon \chi_p(\|s_k\|).
  \eeqn
  Otherwise we have that \req{descent}
  and 
  \beqn{relacc-sk-ok}
\left|\barDT_p^f(x_k,s_k,\omega_k)-\Delta T_p^f(x_k,s_k)\right|
\leq \omega_k \barDT_p^f(x_k,s_k,\omega_k)
\eeqn
are satisfied. Moreover, either $\|s_k\| \geq \mu \epsilon^{\frac{1}{p-q+\beta}}$, or 
\begin{eqnarray}
\phi_{m_k,q}^{\delta_k}(s_k)
&\leq & (1+\kappa_\omega)\,
  \max\left[\frac{\vartheta(1-\kappa_\omega)}{(1+\kappa_\omega)^2}\,\epsilon,
  \frac{\theta \|s_k\|^{p-q+\beta}}{(p-q+\beta)!}\right]\chi_q(\delta_k).
  \label{good-phim}
\end{eqnarray}
  hold.
}

\begin{proof} 
We first note that,  
because of \req{minus-mk} and because Step~2.2
imposes \req{descent}, we have that $\barDT_p^f(x_k,s_k,\omega_k)\geq 0$ at the
end of this step. Let us first consider the case where the calls to the VERIFY
algorithm in Step~2.2  and in Step~2.3 both return ${\tt flag}_s = 1$ or ${\tt
flag}_s = 3$ and note that Step~2.1  yields \req{rvareps-j} with
$T_r=T_r^{f}$, $r=p$ and $\{\zeta_j\}_{j=1}^{r}=\{\varepsilon_{j,i_\varepsilon}\}_{j=1}^p$.
Moreover, $\omega=\omega_k\in(0,1)$ so that  we can use Lemma
\ref{verify-l} to analyse the outcome of the above calls to the VERIFY Algorithm.
If $\|s_k\| = 0$,  Lemma \ref{zerostep} ensures that \req{bar-term-q} holds
for $x = x_{\varepsilon}$ for a radius $\delta \in (0,1)$. Otherwise, we have that
$\barDT_p^f(x_k,s,\omega_k)>0$ because of \req{descent} and, since $s_k$ is then
a global minimizer of $m_k$, that
\beqn{DTDm}
\barDT_p^f(x_k,s,\omega_k) - \frac{\sigma_k}{(p+\beta)!}\|s\|^{p+\beta}
\leq \barDT_p^f(x_k,s_k,\omega_k) - \frac{\sigma_k}{(p+\beta)!}\|s_k\|^{p+\beta}
\eeqn
for all $s$. Thus, if $\|s\|\leq \|s_k\|$, then
$\barDT_p^f(x_k,s,\omega_k)\leq \barDT_p^f(x_k,s_k,\omega_k)$.
This implies that
\[
\globmax_{\mystack{x_k+s \in \calF}{\|s\|\leq\|s_k\|}}\barDT_p^f(x_k,s,\omega_k)
= \barDT_p^f(x_k,s_k,\omega_k).
\]
We may now repeat the proof of Lemma~\ref{acc-S1-l} for the cases  ${\tt flag}_s$
$\in \{1,3\}$, with $q$ replaced by $p$ and $\delta_{k-1}$ replaced by
$\|s_k\|$, and
deduce that \eqref{strong-opt} holds.

Assume now that Algorithm~\ref{step2}  terminates  in Step~2.4.
This means that the VERIFY algorithm
invoked in either Step~2.2 or Step~2.3 terminates with ${\tt flag}_s = 2$, and
we deduce from \req{verif-prop-2} that\req{relacc-sk-ok} holds.

Let us now consider the case  $\|s_k\| < \mu \epsilon^{\frac{1}{p-q+\beta}}$
and note that Lemma \ref{abserr-derm-l}  ensures that  
\req{rvareps-j} is satisfied with $T_r=T_r^{m_k}$, $r=q$ and 
$\{\zeta_j\}_{j=1}^{r}=\{3\varepsilon_{j,i_\varepsilon}\}_{j=1}^q$.
Moreover, the triangle inequality gives
\beqn{triang}
\Delta T_q^{m_k}(s_k,d)
\leq \barDT_q^{m_k}(s_k,d,\omega_k)
    +\left|\barDT_q^{m_k}(s_k,d,\omega_k) - \Delta T_q^{m_k}(s_k,d)\right|.
\eeqn

First, assume that, in Step~2.4, Algorithm~\ref{step2} terminates because
${\tt flag}_d= 1$ is returned by VERIFY. Then,
$\barDT_q^{m_k}(s_k,d_k^{m_k},\omega_k)=0$. Moreover, using \req{triang}, the
definition of $d_k^{m_k}$ given at Step~2.2 of Algorithm \ref{step2},
\req{verif-prop-1} and recalling that $\omega_k\leq 1$,
we obtain that, for all $d$ with $\|d\|\leq\delta_k$, 
\begin{eqnarray}
\Delta T_q^{m_k}(s_k,d)
& \leq & \barDT_q^{m_k}(s_k,d,\omega_k)
         +\left|\barDT_q^{m_k}(s_k,d,\omega_k)-\Delta T_q^{m_k}(s_k,d)\right|\nonumber \\
& \leq & \barDT_q^{m_k}(s_k,d_k^{m_k},\omega_k)
         +\left|\barDT_q^{m_k}(s_k,d,\omega_k)-\Delta T_q^{m_k}(s_k,d)\right|\nonumber\\
&   =  & \left|\barDT_q^{m_k}(s_k,d,\omega_k)- \Delta T_q^{m_k}(s_k,d)\right|\nonumber\\
& \leq & \frac{\vartheta(1-\kappa_\omega)}{2(1+\kappa_\omega)^2}\, \omega_k\,
         \epsilon\,\chi_q(\|d\|)\nonumber\\
& \leq & \frac{\vartheta(1-\kappa_\omega)}{(1+\kappa_\omega)^2}\,
         \epsilon\,\chi_q(\delta_k).\label{S2-term-0}
\end{eqnarray}
If, instead, termination occurs with VERIFY returning ${\tt flag}_d= 2$,
then we will show that for all $d$ with $\|d\|\leq\delta_k$, 
\beqn{S2-term-1}
\Delta T_q^{m_k}(s_k,d)
  \leq (1+\omega_k) \barDT_q^{m_k}(s_k,d_k^{m_k},\omega_k) 
\leq (1+\omega_k) \bigfrac{\theta\|s_k\|^{p-q+\beta}}{(p-q+\beta)!}\,\chi_q(\delta_k).
\eeqn
\noindent
Indeed, from \req{triang}, \req{verif-prop-2}, \req{omega_bound} and the
definition of $d_k^{m_k}$ at Step~2.2 of Algorithm~\ref{step2},
we obtain for all $d$ with $\|d\|\leq\delta_k$
\begin{eqnarray*}
\Delta T_q^{m_k}(s_k,d) &\leq& (1+\omega_k)\barDT_q^{m_k}(s_k,d_k^{m_k},\omega_k),\\
& \le&  (1+\omega_k) \max\Bigg[0,
  \bigglobmax_{\stackrel{x_k+s_k+d\in \calF}{\|d\|\leq\delta_k}}\barDT_q^{m_k}(s_k,d,\omega_k)\Bigg]\\
& = &(1+\omega_k) \barphi_{m_k,q}^{\delta_k}(s_k,\omega_k),
\end{eqnarray*}
in which the equality follows from the definition \req{barphi-def}.  We can then
conclude, using \req{mterm}, that \req{S2-term-1} holds for all
$d$ with $\|d\|\leq\delta_k$.\\

\noindent Finally, if termination occurs instead because VERIFY returns ${\tt flag}_d=
3$, we deduce from the \req{triang} and \req{verif-prop-3} that, for all $d$
with $\|d\|\leq\delta_k$, 
\beqn{S2-term-2}
\Delta T_q^{m_k}(s_k,d)
\leq \frac{\vartheta(1-\kappa_\omega)}{(1+\kappa_\omega)^2}\, \epsilon\,\chi_q(\delta_k)
\eeqn
Observe now that \req{omega_bound}, \req{phi-def} (for $m_k$ at $s_k$) and
each of \req{S2-term-0}, \req{S2-term-1} or \req{S2-term-2} ensures \req{good-phim}.
\end{proof}

Note that \eqref {strong-opt}
 can be viewed as a stronger optimality condition than
\req{term-q} since it implies that the $p$-th (rather than $q$-th with $q\leq
p$) order Taylor expansion of $f$ around $x_k$ is bounded below by a correctly
scaled multiple of $\epsilon$, and in a possibly larger neighbourhood. It
is thus acceptable to terminate the AR$p$DA algorithm with
$x_\epsilon=x_k$ as stated in Step~2.4 of Algorithm~\ref{step2}.

\subsection{The complexity of a single AR$p$DA iteration}\label{compexity-single-s}

The last part of this section is devoted to bounding the evaluation complexity
of a single iteration of the AR$p$DA algorithm.  The count in (approximate)
objective function evaluations is the simplest: these only occur in Step~3 which
requires at most two such evaluations.

Now observe that evaluations of $\{\overline{\nabla_x^j f}\}_{j=1}^p$ possibly
occur in Steps~1.2 and 2.1. However it is important to note that, within these
steps, the \emph{derivatives are evaluated only if the current values of the
  absolute errors are smaller than that used for the previous evalutions} of
the same derivative at the same point ($x_k$). Moreover, these absolute errors
are, by construction, linearly decreasing with rate $\gamma_\varepsilon$
within the same iteration of the AR$p$DA algorithm (they are initialized in
Step~1.1, decreased each time by a factor $\gamma_\varepsilon$ in
\req{decrease-vareps} invoked in Step~1.5, down to values
$\{\varepsilon_{j,i_\varepsilon}\}_{j=1}^p$ which are then passed to the
modified Step~2, and decreased there further in \req{decrease-vareps-2} in
Step~2.5, again by successive multiplication with
$\gamma_\varepsilon$). Furthermore, we have argued already, both for the
modified Step~1 and the modified Step~2, that any of these algorithms
terminates as soon as \req{max-abs-acc} holds for the relevant value of $\xi$,
which we therefore need to determine.  For Step~1, this value is
$\half\omega_k\epsilon$, while, for Step~2, it is given by
\beqn{wrong-xi-2}
\min\left[ \half\omega_k\epsilon,
 \frac{\vartheta(1-\kappa_\omega)}{(1+\kappa_\omega)^2}\frac{\omega_k\epsilon}{2}
 \right]= \frac{\vartheta(1-\kappa_\omega)}{2(1+\kappa_\omega)^2}\,\omega_k\epsilon
\eeqn
when $\|s_k\| < \mu\epsilon^{\frac{1}{p-q+\beta}}$ and by
$\half\omega_k\epsilon$ when $\|s_k\| \geq \mu\epsilon^{\frac{1}{p-q+\beta}}$.
\noindent
As a consequence, we obtain the following lemma.

\llem{one-it-compl-l}{ Suppose that $\omega_k \geq \omega_{\min} > 0$ for all $k$.
  Then each iteration of the AR$p$DA algorithm involves at most 2 (approximate)
  evaluations of the objective function and
  at most $1 + \nu_{\max}(\epsilon)$
  (approximate) evaluations of its $p$ first derivatives, where
  \beqn{nuk-bound}
  \nu_{\max}(\epsilon) = \left\lfloor
  \frac{1}{\log(\gamma_\varepsilon)}
  \left\{\log\left(
  \frac{\vartheta(1-\kappa_\omega)}{6(1+\kappa_\omega)^2}\omega_{\min}\epsilon
  \right)-\log(\kappa_\varepsilon)\right\}
  \right\rfloor.
  \eeqn
}

\begin{proof}
The upper bound on the (approximate) function evaluations immediately
follows from the observation that, as mentioned at the beginning of the
current paragraph, these computations occur at most twice in Step~3 of
Algorithm \ref{algo}. Concerning the second part of the thesis we notice
that, from Lemma \ref{abserr-derm-l}, in Step~2.4 of Algorithm \ref{step2}
we have to make $\{\overline{\nabla_x^j f}(x_k)\}_{j=1}^p$ three times more
accurate than the desired accuracy in $\{\overline{\nabla_d^j
T}{\,}_q^{m_k}(s_k,0)\}_{j=1}^q$, when
$\|s_k\|<\mu\epsilon^{\frac{1}{p-q+\beta}}$ (the input values for the
absolute accuracy values in the VERIFY call are
$\{3\epsilon_j\}_{j=1}^q$). Thus, the VERIFY Algorithm stops whenever
\[
\max_{j\in\{1,\ldots,q\} }\varepsilon_j
\leq \bigfrac{\vartheta(1-\kappa_\omega)\omega_k\epsilon}
             {6(1+\kappa_\omega)^2}.
\]
We may thus conclude from Lemma~\ref{verify-l} that no further reduction in
$\{\varepsilon_j\}_{j=1}^p$ (and hence no further approximation of
$\{\overline{\nabla_x^j f}(x_k)\}_{j=1}^p$) will occur once $i_\varepsilon$,
the number of decreases in $\{\varepsilon_j\}_{j=1}^p$, is large enough to
ensure that
\[
\gamma_\varepsilon^{i_\epsilon} [\max_{j\in\ii{p}}\varepsilon_{j,0}]
\leq  \bigfrac{\vartheta(1-\kappa_\omega)}{6(1+\kappa_\omega)^2}
    \omega_{\min}\epsilon
\]
(Note that this inequality could hold for $i_\epsilon=0$.)
Because of our assumption that
$\omega_k \geq \omega_{\min} $ and \req{upp-vareps}, the above inequality is
then verified when
\[
  i_{\epsilon}\leq\left\lfloor
  \frac{1}{\log(\gamma_\varepsilon)}
  \left\{\log\left(
  \bigfrac{\vartheta(1-\kappa_\omega)}{6(1+\kappa_\omega)^2}\omega_{\min}\epsilon
  \right)-\log(\kappa_\varepsilon)\right\}
  \right\rfloor,
\]
which concludes the proof when taking into account that the derivatives must
be computed at least once per iteration.
\end{proof}\\

\noindent
Note that, for simplicity, we have ignored the fact that only $q \leq p$
derivatives need to be evaluated in Steps~1.2. Lemma~\ref{one-it-compl-l} can
obviously be refined to reflect this observation.

We conclude this section by a comment on what happens whenever exact objective
function and derivatives are used. In that case the (exact) derivatives are
computed only once per iteration of the AR$p$DA algorithm (in Step~1.2 for the
first $q$ and in Step~2.1 for the remaining $p-q$) and every other call to
VERIFY returns ${\tt flag} = 1$ or ${\tt flag}=2$. Moreover, there is no need
to recompute $\barf$ to obtain \req{Df-DT} in Step~3. The evaluation
complexity of a single iteration of the AR$p$DA algorithm then reduces to a
single evaluation of $f$ and its first $p$ derivatives (and
$\nu_{\max}(\epsilon)=1$ for all $k$), as expected.

\numsection{Evaluation complexity of the deterministic AR$p$DA}
\label{complexity-det-s}

This section is devoted to the evaluation complexity analysis of
the AR$p$DA algorithm in the deterministic context. We start by providing a
simple lower bound on the model decrease. 

\llem{Dm-lemma}{\cite[Lemma~3.1]{CartGoulToin18b}
The mechanism of the AR$p$DA algorithm guarantees that, for all $k  \geq 0$,
\beqn{Dphi}
\barDT_p^f(x_k,s_k,\omega_k) > \frac{\sigma_k}{(p+\beta)!} \|s_k\|^{p+\beta},
\eeqn
and so \req{rhokdef2} is well-defined.
}

\begin{proof}We have that
\[
0 < m_k(0) - m_k(s_k)
= \barT_p(x_k,0,\omega_k) - \barT_p(x_k,s_k,\omega_k)
  -\frac{\sigma_k}{(p+\beta)!} \|s_k\|^{p+\beta}.
\]
\end{proof} %epr

\noindent
We next show that the regularization parameter $\sigma_k$ has to remain
bounded, even in the presence of inexact computation of $f$ and its
derivatives. This lemma hinges heavily on \req{barDT-acc}, \req{Df+-DT} and
\req{Df-DT}.

\llem{sigmaupper-lemma}{
Let $f\in \mathcal{C}^{p,\beta}(\Re^n)$.
Then, for all $k\geq 0$,
\beqn{sigmaupper}
\sigma_k
\leq \sigma_{\max} \eqdef \max\left[ \sigma_0,\frac{\gamma_3 (L+3)}{1-\eta_2}\right]
\vspace*{-2mm}
\eeqn
and
\vspace*{-2mm}
\beqn{omegalower}
\omega_k \geq \omega_{\min} \eqdef \min\left[\kappa_\omega,\frac{1}{\sigma_{\max}}\right].
\eeqn
}

\begin{proof}
Assume that
\beqn{siglarge}
\sigma_k \geq \frac{L+3}{1-\eta_2}.
\eeqn
Also observe that, because of the triangle
inequality, \req{relacc-sk-ok} (as ensured by Lemma~\ref{good-step2}) and \req{Df-DT},
\[
\begin{array}{lcl}
|\barT_p^f(x_k,s_k,\omega_k) - T_p^f(x_k,s_k)|
& \leq & |\barf_k(x_k,\omega_k)-f(x_k)| \\*[1ex]
&      &  \hspace*{2cm} + |\barDT_p^f(x_k,s_k,\omega_k)-\Delta T_p^f(x_k,s_k)|\\*[1ex]
& \leq & 2 \omega_k |\barDT_p^f(x_k,s_k,\omega_k)|
\end{array}
\]
and hence, again using the triangle inequality, \req{Df+-DT}, \req{resf},
\req{new-acc}, \req{Dphi} and \req{siglarge},
\[
\begin{array}{ll}
|\rho_k - 1|
& \leq \bigfrac{|\barf_k(x_k+s_k,\omega_k) - \barT_p^f(x_k,s_k,\omega_k)|}
        {\barDT_p^f(x_k,s_k,\omega_k)} \\*[2ex]
& \leq \bigfrac{1}{\barDT_p^f(x_k,s_k,\omega_k)}
   \Big[|\barf_k(x_k+s_k,\omega_k) - f(x_k+s_k)| + |f(x_k+s_k)-T_p^f(x_k,s_k)|\\*[2ex]
&\hspace*{50mm}          + | \barT_p^f(x_k,s_k,\omega_k)-T_p^f(x_k,s_k)| \Big] \\*[2ex]
& \leq \bigfrac{1}{ \barDT_p^f(x_k,s_k,\omega_k)}
   \Big[ |f(x_k+s_k)-T_p^f(x_k,s_k)|+3\omega_k \barDT_p^f(x_k,s_k,\omega_k)\Big]\\*[2ex]
& \leq \bigfrac{1}{\barDT_p^f(x_k,s_k,\omega_k)}
   \Big[ \bigfrac{L}{(p+\beta)!} \|s_k\|^{p+\beta}
   + \bigfrac{3 \barDT_p^f(x_k,s_k,\omega_k)}{\sigma_k} \Big] \\*[2ex]
& < \bigfrac{L}{\sigma_k} + \bigfrac{3}{\sigma_k} \\*[2ex]
& \le 1-\eta_2
\end{array}
\]
and thus that $\rho_k \ge \eta_2$. Then iteration $k$ is very successful in that
$\rho_k \geq \eta_2$ and, because of \req{sigupdate}, $\sigma_{k+1}\leq
\sigma_k$.  As a consequence, the mechanism of the algorithm ensures that
\req{sigmaupper} holds. Observe now that this result and \req{new-acc} imply that,
for all $k$, $\omega_k$ may be chosen such that
$\min[\kappa_\omega,\sigma_{\max}^{-1}]\leq \omega_k \leq \kappa_\omega$,
yielding \req{omegalower}.
\end{proof}%epr

\noindent
It is important to note that \req{omegalower} in this lemma provides the lower
bound on $\omega_k$ required in Lemma~\ref{one-it-compl-l}.
We now borrow a technical result from \cite{CartGoulToin18b}.

\llem{njsp-l-a}{\cite[Lemma~2.4]{CartGoulToin18b}
Let $s$ be a vector of $\Re^n$ and $p\in \Na_0$ and $\beta \in (0,1]$ such that
$j \in \iibe{0}{p}$.  Then
\beqn{npsp-a}
\|\, \nabla_s^j \big(\|s\|^{p+\beta} \big) \, \|_{[j]}
\leq \frac{(p+\beta)!}{(p-j+\beta)!}\|s\|^{p-j+\beta}.
\eeqn
}

\noindent
Our next move is to prove a lower bound on the step norm.  While the proof of
this result is clearly inspired from that of \cite[Lemma~3.3]{CartGoulToin18b},
it nevertheless crucially differs when approximate values are considered
instead of exact ones.

\llem{longs-j-lemma}{
Let $f\in \mathcal{C}^{p,\beta}(\Re^n)$. Then, for all  $k\geq 0$ such that
the AR$p$DA algorithm does not terminate at iteration $k+1$,
\vspace*{-2mm}
\beqn{longs-j}
\|s_k\|\geq \kappa_s \epsilon^{\frac{1}{p-q+\beta}},
\vspace*{-3mm}
\eeqn
where
\vspace*{-1mm}
\beqn{kappas-def}
\kappa_s \eqdef
\min\left\{\mu,\left[\frac{(1-\kappa_\omega)(1-\vartheta)(p-q+\beta)!}
{(1+\kappa_\omega)(L+\sigma_{\max}+\theta(1+\kappa_\omega))}\right]^{\frac{1}{p-q+\beta}}
\right\}.
\eeqn
}

\begin{proof}
If $\|s_k\| \ge \mu \epsilon^{\frac{1}{p-q+\beta}}$, the result is
obvious.  Suppose now that
\beqn{small-s}
\|s_k\| < \mu \epsilon^{\frac{1}{p-q+\beta}}.
\eeqn
Since the algorithm does not terminate at iteration $k+1$, we have that
\vspace*{-2mm}
\[
\barphi_{f,q}^{\delta_k}(x_{k+1}) > \frac{\epsilon}{1+\omega_k} \, \chi_q(\delta_k)
\vspace*{-2mm}
\]
and therefore, using \req{good-phi}, that
\vspace*{-2mm}
\beqn{big-phi-plus}
\phi_{f,q}^{\delta_k}(x_{k+1}) > \frac{1-\omega_k}{1+\omega_k}\, \epsilon\, \chi_q(\delta_k).
\vspace*{-2mm}
\eeqn
Let the global minimum in the definition of $\phi_{f,q}^{\delta_k}(x_{k+1})$ be
achieved at $d$ with $\|d\|\leq \delta_k$.  Then, using \req{phi-def}, the triangle
inequality and \req{npsp-a}, we deduce that 
\begin{eqnarray}
\phi_{f,q}^{\delta_k}(x_{k+1})
&   =  & - \bigsum_{\ell=1}^q \bigfrac{1}{\ell!}\nabla^\ell_xf(x_{k+1})[d]^\ell
         \nonumber\\
& \leq & \left| \bigsum_{\ell=1}^q \bigfrac{1}{\ell!}\nabla^\ell_xf(x_{k+1})[d]^\ell
         - \bigsum_{\ell=1}^q \bigfrac{1}{\ell!}\nabla^\ell_sT^f_p(x_k,s_k)[d]^\ell\right| 
         - \bigsum_{\ell=1}^q \bigfrac{1}{\ell!}\nabla^\ell_sT^f_p(x_k,s_k)[d]^\ell
         \nonumber\\
& \leq & \bigsum_{\ell=1}^q \bigfrac{1}{\ell!}
         \Big[\|\nabla^\ell_xf(x_{k+1})-\nabla^\ell_sT^f_p(x_k,s_k)\|_{[\ell]}\Big]\delta_k^\ell
         \nonumber\\
&      & \hspace*{2cm} - \bigsum_{\ell=1}^q \bigfrac{1}{\ell!}
         \left(\nabla^\ell_s\left[ T^f_p(x_k,s)
             +\frac{\sigma_k}{(p+\beta)!}\|s\|^{p+\beta}\right]_{s=s_k}\right)[d]^\ell
         \nonumber \\
&      & \hspace*{2cm} + \bigsum_{\ell=1}^q
         \bigfrac{\sigma_k}{\ell!(p-\ell+\beta)!}\|s_k\|^{p-\ell+\beta}\delta_k^\ell.
         \label{big-bound1}
\end{eqnarray}         
Now, because of \req{model}, \req{phi-def} (for $m_k$ at $s_k$)
and the fact that $\|d\|\leq \delta_k$, we have that
\[
- \bigsum_{\ell=1}^q \bigfrac{1}{\ell!}\left(\nabla^\ell_s\left[
  T^f_p(x_k,s)+\frac{\sigma_k}{(p+\beta)!}\|s\|^{p+\beta}\right]_{s=s_k}\right)[d]^\ell
= \Delta T_q^{m_k}(s_k,d)\leq \phi_{m_k,q}^{\delta_k}(s_k).
\]
Then, as $\|s_k\|< \mu\epsilon^{\frac{1}{p-q+\beta}}<1$
because of \req{small-s},  we may use \req{good-phim} (ensured by
Lemma~\ref{good-step2}) and \req{resder} and distinguish the
cases where the maximum in \req{good-phim} is attained in its first or its
second argument. In the latter case,  we deduce from \req{big-bound1}  that
\begin{eqnarray}  
\phi_{f,q}^{\delta_k}(x_{k+1})
& \leq & \bigsum_{\ell=1}^q \bigfrac{L}{\ell!(p-\ell+\beta)!}\|s_k\|^{p-\ell+\beta}\delta_k^\ell
         + (1+\kappa_\omega)\bigfrac{\theta\,\chi_q(\delta_k)}{(p-q+\beta)!}
         \|s_k\|^{p-q+\beta}\nonumber\\
&      & \hspace*{2cm} + \bigsum_{\ell=1}^q
         \bigfrac{\sigma_k}{\ell!(p-\ell+\beta)!}\|s_k\|^{p-\ell+\beta}\delta_k^\ell \nonumber\\
& \leq & \bigfrac{\Big[L+\sigma_k+\theta(1+\kappa_\omega)\Big]\chi_q(\delta_k)}
                 {(p-q+\beta)!}\|s_k\|^{p-q+\beta}\label{big-bound};
\end{eqnarray}         
otherwise, \req{big-bound1} guarantees that
\begin{eqnarray}  
\phi_{f,q}^{\delta_k}(x_{k+1})
& \leq & \bigfrac{(L+\sigma_k)\chi_q(\delta_k)}
    {(p-q+\beta)!}\|s_k\|^{p-q+\beta}+\frac{\vartheta(1-\kappa_\omega)}
    {1+\kappa_\omega}\,\epsilon \chi_q(\delta_k).\label{bigb2}
\end{eqnarray}   
Using now \req{big-phi-plus}, \req{omega_bound}, \req{small-s},
\req{big-bound} and \req{bigb2},we thus have that  
\[
\begin{array}{lcl}
\!\!\!\!\|s_k\|
&\!\!\!\!\!\! \geq\!\!\!\!\!\! &\min\!\left\{\!\mu  \epsilon^{\frac{1}{p-q+\beta}},\!
\left[\!\bigfrac{\epsilon(1-\kappa_\omega)(p-q+\beta)!}
           {(1+\kappa_\omega)(L+\sigma_k+\theta(1+\kappa_\omega))}
\!\right]^{\frac{1}{p-q+\beta}}\!\!\!,
\left[\!\bigfrac{\epsilon(1-\kappa_\omega)(1-\vartheta)(p-q+\beta)!}
{(1+\kappa_\omega)(L+\sigma_k)}\!\right]^{\frac{1}{p-q+\beta}} \! \right\}\\
&\!\!\!\!\!\!\geq\!\!\!\!\!\! &\min\!\left\{\!\mu \epsilon^{\frac{1}{p-q+\beta}},
\left[\!\bigfrac{\epsilon(1-\kappa_\omega)(1-\vartheta)(p-q+\beta)!}
           {(1+\kappa_\omega)(L+\sigma_k+\theta(1+\kappa_\omega))}
\!\right]^{\frac{1}{p-q+\beta}}  \!\right\} ,
\end{array}
\]
where we have used the fact that $\theta \in (0,1)$ to obtain the last inequality.
Then \req{longs-j}  follows from \req{sigmaupper}.
\end{proof} %epr

We now combine all the above results to deduce an upper bound on the maximum
number of successful iterations, from which a final complexity bound
immediately follows.

\lthm{upper-theorem}{
Let $f\in \mathcal{C}^{p,\beta}(\Re^n)$ and $\epsilon \in (0,1)$ be given.
Then the AR$p$DA algorithm using the modified Steps~1 
(\hspace*{-1.2mm}\vpageref{step1}) and 2 (\hspace*{-1.2mm}\vpageref{step2})
produces an iterate $x_\epsilon$ such that \req{term-q} or \req{strong-opt}
holds in at most 
\beqn{compl-iters}
\left \lfloor \kappa_p ( f(x_0)- \flow)
\left( \epsilon^{-\frac{p+\beta}{p-q+\beta}} \right)
\right \rfloor+1
\eeqn
successful iterations,
\beqn{ttau}
\tau(\epsilon)
\eqdef \left \lfloor \left\{\left \lfloor
\kappa_p ( f(x_0)- \flow)
\left( \epsilon^{-\frac{p+\beta}{p-q+\beta}} \right)+1
\right \rfloor
\left(1+\frac{|\log\gamma_1|}{\log\gamma_2}\right)+
     \frac{1}{\log\gamma_2}\log\left(\frac{\sigma_{\max}}{\sigma_0}
\right)\right\}\right \rfloor
\eeqn
iterations in total, $2\tau(\epsilon)$ (approximate) evaluations of $f$ and
$(1+\nu_{\max}(\epsilon)) \tau(\epsilon)$
approximate evaluations of $\{\nabla_x^jf\}_{j=1}^p$,
where $\sigma_{\max}$ is given by \req{sigmaupper}, $\omega_{\min}$ by \req{omegalower},
$\nu_{\max}(\epsilon)$ by \req{nuk-bound}, and where
\beqn{kappap-def}
\kappa_p \eqdef \frac{(p+\beta)!}{\eta_1 (1-\alpha)\sigma_{\min}}
\max\left\{
\frac{1}{\mu^{p+\beta}},\left[\frac{(1+\kappa_\omega)(L+\sigma_{\max}+\theta(1+\kappa_\omega))}
         {(1-\kappa_\omega)(1-\vartheta)(p-q+\beta)!}\right]^{\frac{p+\beta}{p-q+\beta}}\right\}.
\eeqn
}

\begin{proof}
At each successful iteration $k$ before termination the algorithm guarantees
the decrease 
\beqn{fdec}
\begin{array}{lcl}
f(x_k)-f(x_{k+1})
& \geq & [\barf_k(x_k,\omega_k) - \barf_k(x_{k+1},\omega_k)]
         - 2 \omega_k\barDT_p^f(x_k,s_k,\omega_k) \\*[2ex]
& \geq & \eta_1 \barDT_p^f(x_k,s_k,\omega_k)
         -\alpha \eta_1\barDT_p^f(x_k,s_k,\omega_k)\\*[2ex]
& > &\bigfrac{\eta_1 (1-\alpha)\sigma_{\min}}{(p+\beta)!} \;\|s_k\|^{p+\beta},
\end{array}
\eeqn
where we used \req{acc-init}, \req{Df+-DT}, \req{Df-DT}, \req{rhokdef2}, \req{Dphi} and
\req{sigupdate}. Moreover  we deduce from \req{fdec} and \req{longs-j} % and \req{sigmaupper}
 that
\beqn{eps1-decr}
f(x_k)-f(x_{k+1}) \geq \kappa_p^{-1} \epsilon^{\frac{p+\beta}{p-q+\beta}}
\tim{where}
\kappa_p^{-1} \eqdef
\bigfrac{\eta_1 (1-\alpha)\sigma_{\min}\kappa_s^{p+\beta} }{(p+\beta)!}.
\eeqn
Thus, since $\{f(x_k)\}$ decreases monotonically,
\[
f(x_0)-f(x_{k+1})
\geq \kappa_p^{-1} \,\epsilon^{\frac{p+\beta}{p-q+\beta}}\,|\calS_k|.
\]
Using that $f$ is bounded below by $f_{\rm low}$, we conclude that
\beqn{Sk1}
| \calS_k |
\leq \kappa_p(f(x_0) - \flow) \epsilon^{-\frac{p+\beta}{p-q+\beta}}
\eeqn
until termination, and the desired bound on the number of successful iterations
follows. Lemma~\ref{SvsU} is then invoked to compute the upper bound on the
total number of iterations, and Lemma~\ref{one-it-compl-l} to deduce the upper
bounds on the number of evaluations of $f$ and its derivatives.
\end{proof}

\noindent
We emphasize that \req{compl-iters} was shown in \cite{CartGoulToin18b} to be
optimal for a quite wide class of minimization algorithms. The slightly
weaker bound $ (1+\nu_{\max}(\epsilon)) \tau(\epsilon)$ may be seen as the
(very modest) price to pay for allowing inexact evaluations.

Focusing on the order in $\epsilon$ and using \req{ttau}, we therefore
obtain the following condensed result on evaluation complexity
for nonconvex optimization.

\lthm{final-in-order}{
Let $f\in \mathcal{C}^{p,\beta}(\Re^n)$. Then, given $\epsilon \in (0,1)$,
the AR$p$DA algorithm
using the modified Steps~1 (\hspace*{-1.2mm}\vpageref{step1}) and 2
(\hspace*{-1.2mm}\vpageref{step2})
needs at most
\vspace*{-2mm}
\[
O\left(\epsilon^{-\frac{p+\beta}{p-q+\beta}}\right)
\tim{iterations and (approximate) evaluations of $f$}
\vspace*{-3mm}
\]
and at most
\vspace*{-2mm}
\[
O\left(|\log(\epsilon)| \epsilon^{-\frac{p+\beta}{p-q+\beta}} \right)
\tim{(approximate) evaluations of the  $p$ first derivatives}
\]
to compute an $(\epsilon,\delta)$-approximate $q$-th-order-necessary minimizer
for the set-constrained problem \req{problem}.
}

\noindent
In particular, if the $p$-th derivative of $f$ is assumed to be globally
Lipschitz rather than merely H\"{o}lder continuous (i.e.\ if $\beta = 1$),
these orders reduce to
\vspace*{-1mm}
\[
O\left(\epsilon^{-\frac{p+1}{p-q+1}}\right)
\tim{iterations and (approximate) evaluations of $f$}
\vspace*{-3mm}
\]
and at most
\vspace*{-1mm}
\[
O\left(|\log(\epsilon)| \epsilon^{-\frac{p+1}{p-q+1}} \right)
\tim{(approximate) evaluations of the $p$ first derivatives.}
\]
As indicated in the comment at the end of Section~\ref{rel-err-s}, all
$O(|\log(\epsilon)|)$ terms reduce to a constant independent of $\epsilon$ if
exact evaluations of $f$ and its derivatives are used, and the above results
then recover the optimal complexity results of \cite{CartGoulToin18b}.

We conclude this section by commenting on the special case where the objective
function evaluations are exact and that of the derivatives inexact. We first
note that this case is already covered by the theory presented above (since
\req{Df+-DT} and \req{Df-DT} automatically holds as their left-hand side is
identically zero), but this remark also shows that the AR$p$DA algorithm can
be simplified by replacing the computation of $\barf(x_k+s_k,\omega_k)$ by
that of $f(x_k+s_k)$ and by skipping the verification and possible
recomputation of $\barf(x_k,\omega_k)$ entirely.  As consequence, the AR$p$DA
algorithm only evaluates the exact objective function $f$ \emph{once} per
iteration, and the maximum number of such evaluations is therefore given by
$\tau(\epsilon)$ instead of $2\tau(\epsilon)$, while the maximum number of
(inexact) derivatives evaluations is still given by
$(1+\nu_{\max}(\epsilon))\tau(\epsilon)$.

\numsection{A variant of the AR$p$DA algorithm}\label{variant-s}

We now describe a variant of the AR$p$DA algorithm for which an even better
complexity can be proved, but at the price of a more restrictive dynamic
accuracy strategy.  In the Step~1.0 of the AR$p$DA algorithm, we allow the
choice of an arbitrary set of $\{\varepsilon_{j,0}\}_{j=1}^p$ with the
constraint that $\varepsilon_{j,0} \leq \kappa_\varepsilon$ for
$j\in\ii{p}$. This allows these accuracy thresholds to vary
non-monotonically from iteration to iteration, providing considerable
flexibility and allowing large inaccuracies even if these thresholds were
made small in past iterations due to local nonlinearity.  A different, more
rigid, strategy is also possible: suppose that the thresholds
$\{\varepsilon_{j,0}\}_{j=1}^p$ are not reset at each iteration, that is
\beqn{variant-def}
  \tim{Step~1.1 is only executed for $k=0$.}
\eeqn
This clearly results in a monotonic decrease of each $\varepsilon_j$ across
all iterations. As a consequence, $\nu_{\max}(\epsilon)$ in \req{nuk-bound}
now bounds the total number of reductions of the $\varepsilon_j$ over all
iterations, i.e. on the \emph{entire run of the algorithm}.
We then deduce that the total number of derivatives evaluation is then bounded by
$\nu_{\max}(\epsilon) + \tau(\epsilon)$ (instead of
$(1+\nu_{\max}(\epsilon))\tau(\epsilon)$) and we may establish the
worst-case complexity of the resulting ``monotonic'' variant as follows.

\lthm{upper-theorem-variant}{ Let $f\in \mathcal{C}^{p,\beta}(\Re^n)$ and
  $\epsilon \in (0,1)$ be given.  Then the AR$p$DA algorithm using the
  modified Steps~1 (\hspace*{-1.2mm}\vpageref{step1}) and 2
  (\hspace*{-1.2mm}\vpageref{step2}) as well as the modified rule
  \req{variant-def} produces an iterate $x_\epsilon$ such that \req{term-q} or
  \req{strong-opt} holds in at most \req{compl-iters} successful iterations,
  $\tau(\epsilon)$ iterations in total, $2\tau(\epsilon)$ (approximate)
  evaluations of $f$ and $\nu_{\max}(\epsilon) + \tau(\epsilon)$ approximate
  evaluations of $\{\nabla_x^jf\}_{j=1}^p$, where $\tau(\epsilon)$ is given by
  \req{ttau}, $\kappa_p$ is given by \req{kappap-def}, $\sigma_{\max}$ by
  \req{sigmaupper}, $\omega_{\min}$ by \req{omegalower} and
  $\nu_{\max}(\epsilon)$ by \req{nuk-bound}.
}

As above, this complexity bound can be condensed to 
\beqn{more-final-in-order}
\begin{array}{l}
  O\left(\epsilon^{-\frac{p+\beta}{p-q+\beta}}\right)
  \tim{iterations and (approximate) evaluations of $f$}\\
  O\left(|\log(\epsilon)| + \epsilon^{-\frac{p+\beta}{p-q+\beta}} \right)
  \tim{(approximate) evaluations of the $p$ first derivatives,}
\end{array}
\eeqn
typically improving on Theorem~\ref{final-in-order}.
When $p=2$, $q=1$ and $\beta=1$, the AR$p$DA variant using
the more restrictive accuracy strategy \req{variant-def} requires at most 
\[
O\Big( |\log(\epsilon)| + \epsilon^{-3/2} \Big)
\]
(approximate) evaluations of the gradient, which corresponds to the bound
derived for the ARC-DFO algorithm of \cite{CartGoulToin12a}.  This is not
surprising as this latter algorithm uses a monotonically decreasing sequence
of finite-difference stepsizes, implying monotonically decreasing
gradient-accuracy thresholds.

One should however notice that the improved bound \req{more-final-in-order}
comes at the price of asking, for potentially many iterations, an accuracy on
$\{\nabla_x^jf\}_{j=1}^p$ which is tighter than what is needed to ensure
progress of the minimization.  In practice, this might be a significant
drawback. We will thus restrict our attention, in what follows, to the
original AR$p$DA algorithm, but similar developments are obviously possible
for the ``monotonic'' variant just discussed.

\numsection{Application to unconstrained and bound-constrained first- and
  second-order nonconvex inexact optimization}
\label{first-order-s}

Because of its wide-ranging applicability, the framework discussed above may
appear somewhat daunting in its generality.  Moreover, the fact that it
involves (possibly constrained) global optimization subproblems in several of
its steps may suggest that it has to remain conceptual. We show in this
section that \emph{this is not the case}, and stress that it is \emph{much
simpler when  specialized to small values of $p$ and $q$} (which are, for
now, the most practical ones) and that our approach leads to elegant and
implementable numerical algorithms. To illustrate this point, we now review
what happens for $p\leq 2$.

We first discuss the case where one seeks to compute a first-order
critical point for an unconstrained optimization problem using approximate
function values as well as approximate first derivatives.  For simplicity of
exposition, we will also assume that the gradient of $f$ is Lipschitz (rather
than H\"{o}lder) continuous. In our general context, this means that we
consider the case where $q =1$, $p = 1$, $\beta = 1$ and $\calF = \Re^n$.
We first note that, as pointed out in \req{phi1-def},
\beqn{phi-barphi-1rst}
\phi_{f,1}^\delta(x) = \|\nabla_x^1f(x)\|\delta
\tim{ and }
\barphi_{f,1}^\delta = \|\overline{\nabla_x^1f}(x)\|\delta
\tim{ irrespective of } \delta \in (0,1],
  \eeqn
which means that, since we can choose $\delta=1$, Step~1 of
the AR$p$dA algorithm reduces to the computation of an approximate gradient
$\overline{\nabla_x^1f}(x_k)$ with relative error $\omega_k$ and verification
that $\epsilon$-approximate optimality is not yet achieved. If that is the
case, computing $s_k$ at Step~2 is also extremely simple since it is easy to
verify that
\[
s_k = s_k^* = -\frac{1}{\sigma_k}\overline{\nabla_x^1f}(x_k).
\]
Lemma~\ref{step-ok-l} then ensures that this step is acceptable for some
$\delta_k \in (0,1]$, the value of which being irrelevant since it is not used
in Step~1 of the next iteration.  Moreover, if the relative error on
$\overline{\nabla_x^1f}(x_k)$ is bounded by $\omega_k$, then
\[
\begin{array}{lcl}
|\barDT_1^f(x_k,s_k)- \Delta T_1^f(x_k,s_k)|
& \leq & \|\overline{\nabla_x^1f}(x_k)- \nabla_x^1f(x_k)\|
         \bigfrac{\|\overline{\nabla_x^1f}(x_k)\|}{\sigma_k}\\*[2ex]
& \leq  &\omega_k \bigfrac{\|\overline{\nabla_x^1f}(x_k)\|^2}{\sigma_k} \\*[2ex]
&   =   & \omega_k \barDT_1^f(x_k,s_k)
\end{array}
\]
and \req{barDT-acc} automatically holds, so that no iteration is needed in
Algorithm~\ref{step2}.  The resulting algorithm, where we have made the
modified Step~1  explicit, is given as Algorithm~\ref{ARpDA1} (AR$1$DA)
\vpageref{ARpDA1}.
\newpage

\algo{ARpDA1}{The AR1DA Algorithm}
{
\vspace*{-0.3 cm}
\begin{description}
\item[Step 0: Initialization.]
  An initial point $x_0\in\Re^n$ and an initial regularization parameter
  $\sigma_0>0$ are given, as well as an accuracy level $\epsilon \in (0,1)$
  and an initial relative accuracy $\omega_0\geq 0$.  The constants  $\alpha$,
  $\kappa_\omega$, $\kappa_\varepsilon$,  $\eta_1$, $\eta_2$, $\gamma_1$, $\gamma_2$, 
  $\gamma_3$ and $\sigma_{\min}$ are also given and satisfy
  $\sigma_{\min} \in (0, \sigma_0]$,
\[
0 < \eta_1 \leq \eta_2 < 1, \;\;
0< \gamma_1 < 1 < \gamma_2 < \gamma_3,
\]
\vspace*{-8mm}
\[
\kappa_\varepsilon \in (0,1]
\ms
\alpha \in (0,1),
\ms
\kappa_\omega \in (0,\half \alpha \eta_1]
\tim{ and }
\omega_0 =\min\left[\kappa_\omega,\frac{1}{\sigma_0}\right].
\]
Set $k=0$.

\item[Step 1: Compute the optimality measure and check for termination. ]
  Initialize $\varepsilon_{1,0} = \kappa_\varepsilon$ and set $i=0$.
  Do
  \begin{enumerate}
    \item compute $\overline{\nabla_x^1f}(x_k)$ with
      $\|\overline{\nabla_x^1f}(x_k)-\nabla_x^1f(x_k)\| \leq
      \varepsilon_{1,i}$ and increment $i$ by one.
    \item if $\|\overline{\nabla_x^1f}(x_k)\|\leq \epsilon/(2(1+\omega_k))$,
      terminate with $x_\epsilon=x_k$;
    \item if $\varepsilon_{1,i}\leq \omega_k \|\overline{\nabla_x^1f}(x_k)\|$,
      go to Step 2;
    \item set $\varepsilon_{1,i+1}= \gamma_\varepsilon\varepsilon_{1,i}$ and
      return to item 1 in this enumeration.
    \end{enumerate}
    
\item[Step 2: Step calculation. ]
  Set
  \vspace*{-2mm}
  \[
  s_k = -\overline{\nabla_x^1f}(x_k)/\sigma_k
  \tim{ and }
  \barDT_1^f(x_k,s_k,\omega_k)= \|\overline{\nabla_x^1f}(x_k)\|^2/\sigma_k.
  \vspace*{-2mm}
  \]

\item[Step 3: Acceptance of the trial point. ] \ \\
  Compute $\barf_k(x_k+s_k,\omega_k)$ ensuring that
  \vspace*{-2mm}
  \beqn{Df+-DT-first}
  |\barf_k(x_k+s_k,\omega_k)-f(x_k+s_k)| \leq \omega_k |\barDT_1^f(x_k,s_k,\omega_k)|.
  \eeqn
  Also ensure (by setting $\barf_k(x_k,\omega_k) =
  \barf_{k-1}(x_k,\omega_{k-1})$ or by (re)computing $\barf_k(x_k,\omega_k)$) 
  that 
  \beqn{Df-DT-first}
  |\barf_k(x_k,\omega_k)-f(x_k)| \leq \omega_k |\barDT_1^f(x_k,s_k,\omega_k)|
  \eeqn
  Then define
  \vspace*{-2mm}
  \beqn{rhokdef-first}
  \rho_k=\frac{\barf_k(x_k,\omega_k)-\barf_k(x_k+s_k,\omega_k)}
              {\barDT_1^f(x_k,s_k,\omega_k)}.
  \eeqn
If $\rho_k \geq \eta_1$, then define
$x_{k+1} = x_k + s_k$; otherwise define $x_{k+1} = x_k$.

\item[Step 4: Regularization parameter update. ]
Set
\beqn{sigupdate-first}
\sigma_{k+1} \in \left\{ \begin{array}{ll}
{}[\max(\sigma_{\min}, \gamma_1\sigma_k), \sigma_k ]  & \tim{if} \rho_k \geq \eta_2, \\
{}[\sigma_k, \gamma_2 \sigma_k ]          &\tim{if} \rho_k \in [\eta_1,\eta_2),\\
{}[\gamma_2 \sigma_k, \gamma_3 \sigma_k ] & \tim{if} \rho_k < \eta_1.
  \end{array} \right.
\eeqn

\item[Step 5: Relative accuracy update. ]
Set
\vspace*{-2mm}
\beqn{new-acc-first}
\omega_{k+1} = \min \left[ \kappa_\omega,\frac{1}{\sigma_{k+1}}\right].
\eeqn
%\vspace*{-2mm}
Increment $k$ by one and go to Step~1.
\end{description}
}

\noindent
Theorem~\ref{final-in-order} then guarantees that the AR1DA Algorithm will find an
$\epsilon$-approximate first-order minimizer for the unconstrained version of
problem \req{problem} in at most $O\big(\epsilon^{-2}\big)$ iterations and
approximate evaluations of the objective function
(which is proved in \cite{CartGoulToin18b} to be optimal) and at most
$O\big(|\log(\epsilon)|\epsilon^{-2}\big)$ approximate evaluations of the
gradient. Note that\vskip 5pt
\begin{enumerate}
  \item the accuracy requirement is truly adaptive and the absolute accuracy
    $\varepsilon_{1,i}$ may remain quite large as long as
    $\|\overline{\nabla_x^1f}(x_k)\|$ itself remains large, as shown by item~3
    in Step~1. \vskip 5pt
  \item The accuracy requirement for computing $\barf$ does not depend on the
    absolute accuracy of the gradient, but only on its norm (squared). At
    initial iterations, this may be quite large.\vskip 5pt
  \item The AR1DA Algorithm is very close in spirit
    to the trust-region with dynamic accuracy of \cite[Sections~8.4.1.1 and
    10.6]{ConnGoulToin00} and, when values of $f$ are computed exactly,
    essentially recovers a proposal in \cite{Bonn18}. It is also close to the
    proposal of \cite{PaquSche18}, which is based on an Armijo-like linesearch
    and has similar accuracy requirements. \vskip 5pt
\end{enumerate}

We now turn to the case where one seeks a first-order critical point for an
unconstrained problem using approximate gradients and Hessians (under the
assumption that the exact Hessian is Lispchitz continuous).  As for the case
$p=q=1$, we have that \req{phi-barphi-1rst} holds, making the verification of
optimality in Step~1 relatively easy.  Computing $s_k$ is now more complicated
but still practical, as it now implies minimizing the regularized quadratic
model $m_k$ starting from $x_k$ until a step $s_k$ is found such that
\[
\|s_k\|\geq \mu \epsilon^{\half}
\tim{ or }
\barphi_{m_k,1}^\delta(s_k,\omega_k)
= \|\nabla_s^1 m_k(s_k)\| \leq \half \theta \|s_k\|^2
\]
(as proposed in \cite{CartGoulToin11d}, see also
\cite{Grie81,NestPoly06,CartGoulToin09a,Duss15}), with the additional constraint that, 
for $s_k\neq 0$,
\beqn{abserr-2-1}
\max[\varepsilon_{1,i},\varepsilon_{2,i}]
\leq  \omega_k\frac{\barDT_2^f(x_k,s_k,\omega_k)}{\chi_2(\|s_k\|)}
\eeqn
where
\[
\barDT_2^f(x_k,s_k,\omega_k)
= -\overline{\nabla_x^1f}(x_k)^Ts_k- \half s_k^T \overline{\nabla_x^2f}(x_k)s_k.
\]

The resulting algorithm AR2DA is quite
similar to AR1DA and is omitted for brevity.  We note that\vskip 5pt
\begin{enumerate}
\item Algorithm AR2DA is guaranteed by Theorem~\ref{final-in-order} to find an
  $\epsilon$-approximate first-order minimizer for the unconstrained version
  of problem \req{problem} in at most $O\big(\epsilon^{-3/2}\big)$ iterations
  and approximate evaluations of the objective function (which is proved in
  \cite{CartGoulToin18b} to be optimal) and at most
  $O\big(|\log(\epsilon)|\epsilon^{-3/2}\big)$ approximate evaluations of the
  gradient and Hessian. \vskip 5pt 
\item As for AR1DA, the absolute accuracies required by AR2DA on the
  approximate function, gradient and Hessian only depend on the magnitude
  of the Taylor increment, which is typically quite large in early
  iterations. The relative errors on the latter two remain bounded away from
  zero.\vskip 5pt
\item The absolute accuracies required on the approximate gradient and Hessian
  are comparable in magnitude, although \req{abserr-2-1} could be exploited to
  favour one with respect to the other.\vskip 5pt
\end{enumerate}

The case where $p=2$ and $q=2$ (i.e.\ when second-order solutions are sought)
is also computationally quite accessible: calculating the optimality measure
$\barphi_{f,1}^{\delta_k}(x_k,\omega_k)$ or
$\barphi_{m_k,1}^{\delta_k}(s_k,\omega_k)$ now involve a standard trust-region
subproblem, for which both exact and approximate numerical solvers are known
(see \cite[Chapter~7]{ConnGoulToin00} for instance), but the rest of the
algorithm --- in particular its adaptive accuracy requirement --- is very
similar to what we just discussed (see also
\cite{CartGoulToin18b}). Theorem~\ref{final-in-order} then ensures that
resulting method converges to an $\epsilon$-approximate second-order-necessary
minimizer for the unconstrained version of problem \req{problem} in at most
$O\big(\epsilon^{-3}\big)$ iterations and approximate evaluations of the
objective function and at most $O\big(|\log(\epsilon)|\epsilon^{-3}\big)$
approximate evaluations of the gradient and Hessian.

We conclude this section by a brief discussion of the case where $q=1$ and
$p \in \{1,2\}$ as before, but where $\calF$ is now defined by bound
constraints.  It is clear that evaluating and enforcing such constraints (by
projection, say) has negliglible cost and therefore falls in our framework. In
this case, the calculations of
$\barphi_{f,1}^{\delta_k}(x_k,\omega_k)$ or
$\barphi_{m_k,1}^{\delta_k}(s_k,\omega_k)$ now involve simple linear
optimization problems\footnote{Formerly known as linear programming problems,
  or LPs.}, which is computationally quite tractable. If $p=1$, Step~2.2 and
2.3 involve convex quadratic optimization, while they involve minimizing
a regularized quadratic model if $p=2$. All results remain the same, and
the AR$p$DA algorithm is then guaranteed to find a bound-constrained approximate
first-order approximate minimizer in at most $O\big(\epsilon^{-2}\big)$ or
$O\big(\epsilon^{-3/2}\big)$ iterations and approximate evaluations of the
objective function (which is proved in \cite{CartGoulToin18b} to be optimal)
and at most $O\big(|\log(\epsilon)|\epsilon^{-2}\big)$ or
$O\big(|\log(\epsilon)|\epsilon^{-3/2}\big)$ approximate evaluations of the
gradient and Hessian. The same algorithms and results
obviously extend to the case where $\calF$ is a convex polyhedral set or any
closed non-empty convex set, provided the cost of the projection on this set
remains negligible compared to that of (approximately) evaluating the
objective function and its derivatives.

\numsection{A stochastic viewpoint on AR$p$DA}
\label{complexity-subs-s}
\subsection{Probabilistic complexity}

In this section we consider the case where the bounds 
 $\{\varepsilon_j\}_{j=1}^p$ on the absolute errors on the
derivative tensors $\{\nabla_x^jf(x)\}_{j=1}^p$ are satisfied with
probability at least $(1-t)$, with $t\in (0,1)$.
This may occur, for instance, if the approximate derivative tensors are
obtained by some stochastic sampling scheme, as we detail below.
We therefore assume  that
\beqn{vareps-j-prob}
Pr\Big[\|\overline{\nabla_x^jf}(x_k)-\nabla_x^jf(x_k)\|_{[j]}\leq\varepsilon_j\Big]
\geq (1-t)
\tim{ for each } j \in \iibe{1}{p}.
\eeqn
We also assume that inequalities \eqref{Df+-DT} and \eqref{Df-DT} in Step 3 of
the AR$p$DA algorithm are satisfied with probability at least $(1-t)$, i.e. 
\beqn{Df+-DT-prob}
  Pr\Big[|\barf_k(x_k+s_k,\omega_k)-f(x_k+s_k)| \leq \varepsilon_0\Big]\ge 1-t,
\eeqn
and
\beqn{Df-DT-prob}
Pr\Big[|\barf_k(x_k,\omega_k)-f(x_k)| \leq \varepsilon_0\Big] \ge 1-t
  \eeqn
where we have defined $\varepsilon_0\eqdef \omega_k|\barDT_p^f(x_k,s_k,\omega_k)|$.
Clearly, different values for $t$ could be chosen in
\eqref{vareps-j-prob}, one for each index (tensor order)
$j\in\{1,...,p\}$. Similarly, different values of $t$ in \eqref{Df+-DT-prob}
and \eqref{Df-DT-prob} could be considered. However, for the sake of
simplicity, we assume here that all the inequalities involved in
\eqref{vareps-j-prob}--\eqref{Df-DT-prob} hold with the same fixed lower bound
$(1-t)$ on the probability of success.  We also assume that the events in
\eqref{vareps-j-prob}--\eqref{Df-DT-prob} are independent.

Stochastic variants of trust-region and adaptive cubic regularization
  methods have been ana\-lyzed in \cite{BandScheVice14,BlanCartMeniSche16,
    CartSche17, XuRoosMaho17, YaoXuRoosMaho18}.   In \cite{BlanCartMeniSche16,
    CartSche17}, complexity results are given in expectation, while the
    analysis is carried out in probability in
    \cite{BandScheVice14,XuRoosMaho17, YaoXuRoosMaho18}.
  We choose to follow the high-probability approach of
  \cite{XuRoosMaho17, YaoXuRoosMaho18}, where an overall and cumulative
  success of  \eqref{vareps-j-prob}--\eqref{Df-DT-prob} is assumed along all the
  iterations up to termination.

We stress that Algorithms~\ref{step1} and \ref{step2} terminates independently
of the satisfaction of the accuracy requirements on the tensor
derivatives. This is due to the fact that termination relies on the inequality
\req{max-abs-acc}.  Moreover, during the iterations of either of these
algorithms before the last, it may happen that the accuracy on the tensor
derivatives fails to be achieved, but this has no impact on the worst-case
complexity.  Satisfying the accuracy requirement is only crucial in the last
iteration of Algorithm~\ref{step1} or \ref{step2} (that is in Steps~1.2 and
2.2).  Let $\calE_r(S)$ be the event: ``the relations
\[
\|\overline{\nabla_x^jf}(x_k)-\nabla_x^jf(x_k)\|_{[j]}\leq \varepsilon_j
\tim{for all} j\in\ii{r}
\]
hold for some $j$ at Step~$S$ of the last iteration of the relevant algorithm''.
In Step~1.2, inexact values are computed for the first $q$ derivatives, 
and the probability that event  $\calE_q(1.2)$ occurs 
is therefore at least $(1-t)^q$. Similarly, the probability that event
$\calE_q(2.2)$ occurs is at least $(1-t)^p$. Finally, at Step 3 of the AR$p$DA
algorithm, the probability that both \eqref{Df+-DT} and \eqref{Df-DT} hold is
at least $(1-t)^2$.  Then, letting for $i\in\ii{k}$, $\calE_{[i]}$ be the event:
``Inequalities \req{barDT-acc}, \req{Df+-DT} and \req{Df-DT} hold at iteration
$i$, of the AR$p$DA algorithm'', the probability that 
$\calE_{[i]}$ occurs is then at least $(1-t)^{p+q+2}$.  Finally,
letting $\calE(k)$ be the event: ``$\calE_{[i]}$ occurs for all
iterations $i \in \ii{k}$ of the AR$p$DA algorithm'', we deduce that
\[
Pr\Big[\calE(k)\Big]
\equiv Pr\left[\bigcap_{i=1}^k \calE_{[i]}\right]\geq(1-t)^{k(p+q+2)}.
\]
Thus, requiring  that the event $\calE(k)$ occurs with  probability
at least $1-\bart$, we obtain that
\[
Pr\Big[\calE(k)\Big]\geq(1-t)^{k(p+q+2)} = 1-\bart,
\tim{i.e.,}
t = 1-(1-\bart)^{\frac{1}{k(p+q+2)}}
= O\left(\frac{\bart}{k(p+q+2)}\right).
\]
Taking into account that, when \req{barDT-acc}, \req{Df+-DT} and \req{Df-DT}
hold, the AR$p$DA algorithm terminates in at most
$k=O\bigl(\epsilon^{-\frac{p+\beta}{p-q+\beta}}\bigr)$ iterations  
(as stated by Theorem~\ref{final-in-order}), we deduce the following result.

\lthm{stoch-final}{Let $f\in \mathcal{C}^{p,\beta}(\Re^n)$. Suppose that
the probabilistic assumptions of this section hold and that, at each of
iteration of the AR$p$DA algorithm, the probability  $t$ satisfies 
\begin{equation}
\label{t}
t = O\left(\frac{\bart\,\epsilon^{\frac{p+\beta}{p-q+\beta}}}{(p+q+2)}\right).
\end{equation}
Then, given $\epsilon \in (0,1)$, the conclusions of
Theorem~\ref{final-in-order} hold with probability at least $(1-\bart)$.
}

\noindent
As a consequence, when $p=q=2$ and $\beta=1$ we have to choose
$t=O\left(\sfrac{1}{6} \,\bart\, \epsilon^3\right)$, while, when
$p=q=\beta=1$, we have to choose
$t=O\left(\sfrac{1}{4}\,\bart \, \epsilon^2\right)$.

We stress that the above analysis is unduly pessimistic in the case where
$p=q=1$. Indeed, as already noticed in Section~\ref{first-order-s}, no
reduction in $\{\varepsilon_j\}$ is necessary at Step~2, as \req{barDT-acc} is
automatically enforced whenever the relative error on the first derivative
$\overline{\nabla_x^1 f}(x)$ is bounded by $\omega_k$.  Noting that this last
event has probability at least $1-t$, we can conclude that $Pr(\calE_{[i]})\geq(1-t)^3$
and to get the optimal complexity $O\bigl(\epsilon^{-2}\bigr)$ with
probability at least $1-\bart$, we need to choose $t=O\left(\third\,\bart\,
  \epsilon^{2}\right)$.  We also emphasize that the purpose of
Theorem~\ref{stoch-final} is limited to offer guidance on desirable value of
$t$ and not to prescribe an algorithmically binding bound. Indeed some of the
constants involved in the bound of Theorem~\ref{final-in-order} (and thus of
Theorem~\ref{stoch-final}) are typically unknown a priori (which is why we
have not been more specific in \req{t}).

\subsection{Sample size in subsampling for finite-sum problems}

In what follows, we now focus on the solution of 
large-scale instances of the finite-sum problems
arising in machine learning and data analysis, that are modelled as
\beqn{finite-sum}
\min_{x\in\cal F} \, f(x)=\frac{1}{N}\sum_{i=1}^N{\psi_i(x)},
\eeqn
with $N>0$ and $\psi_i:{\Re}^n\rightarrow {\Re}$.  
Restricting ourselves to the cases where $p\leq 2$, we discuss the application
of Algorithm AR1DA and AR2DA to problem (\ref{finite-sum}). In this case, the
approximation of the objective function's value and of first and second
derivatives is obtained by a subsampling procedures, i.e.\ these quantities
are approximated  by randomly sampling component functions $\psi_i$. More
precisely, at iteration $k$ these approximations take the form: 
\[
\barf_k (x_k,\omega_k)= \frac{1}{|{\cal D}_{k,1}|} \sum_{i \in {\cal D}_{k,1}} \psi_i(x_k),
\ms
\barf_k (x_k+s_k,\omega_k)= \frac{1}{|{\cal D}_{k,2}|} \sum_{i \in {\cal D}_{k,2}} \psi_i(x_k+s_k),
\]
\[
\overline{\nabla_x^1 f}(x_k)
= \frac{1}{|{\cal G}_k|} \sum_{i \in {\cal G}_k} \overline{\nabla_x^1 \psi_i}(x_k), 
\tim{ and }
\overline{\nabla_x^2 f}(x_k)
= \frac{1}{|{\cal H}_k|} \sum_{i \in {\cal H}_k} \overline{\nabla_x^2 \psi_i}(x_k), 
\]
where ${\cal D}_{k,1}$, ${\cal D}_{k,2}$, ${\cal G}_k$ and ${\cal H}_k$ are subsets of
$\{1,2,\ldots,N\}$.  The question then arises of estimating the cardinality of
these sample sets in order to ensure that the approximations of the objective
function's value and its first and second derivatives satisfy
\eqref{vareps-j-prob} for $j=1$ and $j=2$, \req{Df+-DT-prob} and
\req{Df-DT-prob}. This issue can be addressed using the
  operator-Bernstein inequality given in \cite{Trop15} and recently extended
  in \cite{LuoQiToint19} to random tensors of general order. In the next
  theorem we derive our final result concerning the sample sizes for
  subsampling the objective function and its derivatives up to order two.

\lthm{sampling-th}{Suppose that there exist non-negative constants
  $\{\kappa_{\psi,j}\}_{j=0}^2$ such that, for $x\in \Re^n$ and all $j\in\{0,1,2\}$  
  \beqn{der-norm-bound}
  \max_{i \in\ii{N}}\|\nabla_x^j\psi_i(x)\| \leq \kappa_{\psi,j}(x).
  \eeqn
  Let  $t \in (0,1)$ and
  suppose that a subsample $\calA_k$ is chosen randomly and uniformly from
  $\ii{N}$ and that, for some $j\in\{0,1,2\}$, one computes
  \[
  \overline{\nabla_x^j f}(x)
  = \frac{1}{|\calA_k|} \sum_{i \in \calA_k} \overline{\nabla_x^j \psi_i}(x),
  \]
  with
  \beqn{sample-size}
  |\calA_k|
  \geq \min\left \{ N,\left\lceil\frac{4\kappa_{\psi,j}(x)}{\varepsilon_j}
  \left(\frac{2\kappa_{\psi,j}(x)}{\varepsilon_j}+\frac{1}{3}\right)
  \,\log\left(\frac{d}{t}\right)\right\rceil\right \},
  \eeqn
  where
  \[
  d=\left\{\begin{array}{ll}
         2,   & \tim{if} j=0,\\
         n+1, & \tim{if} j=1,\\
         2n,  & \tim{if} j=2.
       \end{array}\right.
  \]

  Then condition \eqref{vareps-j-prob} holds for $x=x_k$ with
  probability at least $(1-t)$ if  $j\in\{1,2\}$, or, if $j=0$, each of the conditions
  \req{Df+-DT-prob} and \req{Df-DT-prob} holds with probability at least
  $(1-t)$ for $x=x_k+s_k$ and $x=x_k$, respectively.
}

\begin{proof}
Let $\calA_k\subseteq \ii{N}$ be a sample set of cardinality
$|\calA_k|$. Consider $j\in \{0,1,2\}$ and $|\calA_k|$ random tensors
$Z_u(x)$ such that, 
\[
Pr\left[Z_u(x)= \nabla_x^j \psi_i(x)\right]
=\frac{1}{N},\quad (i\in\ii{N}).
\]
For $u\in\calA_k$, let us define
\[
X_u \eqdef\left(Z_u(x)-\nabla_x^j f(x)\right),
\qquad
\overline{\nabla_x^j f}(x)
\eqdef \frac{1}{|\calA_k|}\sum_{u\in\calA_k}Z_u(x)
\]
and
\[
X \eqdef\sum_{u\in\calA_k}X_u
=|\calA_k|\left(\overline{\nabla_x^j f}(x)-\nabla_x^j f(x)\right).
\]
Since \req{finite-sum} gives that
\[
\frac{1}{N}\sum_{i=1}^N \nabla_x^j \psi_i(x) =\nabla_x^j f(x),
\]
we deduce that
\[
E(X_u)
=\frac{1}{N}\sum_{i=1}^N\left( \nabla_x^j \psi_i(x)-\nabla_x^j f(x)\right)
=0,
\quad u\in\calA_k.
\]
Moreover, assuming $Z_u(x)=\nabla_x^j \psi_l(x)$ for some $l\in\ii{N}$
and using \req{der-norm-bound}, we have that
\[
\|X_u\|
\leq \Bigg\|\frac{N-1}{N}\,\nabla_x^j \psi_l(x)-\frac{1}{N}\sum_{i\in\ii{N}\setminus\{l\}}\nabla_x^j \psi_i(x)\Bigg\|
\leq 2\,\frac{N-1}{N}\, \kappa_{\psi,j}(x)
\leq 2\kappa_{\psi,j}(x),
\]
so that the variance of $X$ can be bounded as follows:
\[
\begin{split}
v(X)
&=\max \left[\|E(XX^T)\|,\|E(X^TX)\|\right]\\
&=\max \Big[\Big\|\sum_{u\in\calA_k}E(X_uX_u^T)\Big\|,\Big\|\sum_{u\in\calA_k}E(X_u^TX_u)\Big\|\Big]\\
&\leq \max \Big[\sum_{u\in\calA_k}\|E(X_uX_u^T)\|,\sum_{u\in\calA_k}\|E(X_u^TX_u)\|\Big]\\
&\leq \max \Big[\sum_{u\in\calA_k}E(\|X_uX_u^T\|),\sum_{u\in\calA_k}E(\|X_u^TX_u\|)\Big]\\
&\leq \sum_{u\in\calA_k}E(\|X_u\|^2)\le 4 |\calA_k|\kappa_{\psi,j}^2(x),
\end{split}
\]
in which the first and the third inequalities hold because of the triangular
inequality, while the second is due to the Jensen's inequality (note that the
spectral norm $\|\cdot\|$ is convex). Therefore, according to the
Operator-Bernstein Inequality stated in \cite[Theorem~6.1.1]{Trop15}, we
obtain that
\beqn{uprob}
Pr\Big[\|\overline{\nabla_x^j f}(x)-\nabla_x^j f(x)\|\ge \epsilon_j\Big]
= Pr\Big[\|X\|\ge \epsilon_j|\calA_k|\Big]
\leq d\,e^{-\frac{\epsilon_j^2|\calA_k|}{4\kappa_{\psi,j}(x)\left(2\kappa_{\psi,j}(x)+\third\epsilon_j\right)}},
\eeqn
with $d=2$ if $j=0$, $d=n+1$ if $j=1$ and $d=2n$ if $j=2$. Then, bounding the
right-hand side of \req{uprob} by $t$, taking logarithms and extracting
$|\calA_k|$ gives \req{sample-size}.
\end{proof}\\%epr

\noindent
In particular, Theorem~\ref{sampling-th} gives the lower bounds 
\beqn{Dkl}
 |\calD_{k,\ell}|\geq\min \left \{N, \left\lceil
 \frac{4\kappa_{\psi,j}(x)}{\varepsilon_0}
  \left(\frac{2\kappa_{\psi,j}(x)}{\varepsilon_0}+\frac{1}{3}\right)
  \log\left(\frac{2}{t}\right)\right \rceil\right \}, \quad \ell=1,2,
\eeqn
\beqn{Gk}
|\calG_k|\geq\min \left \{N,\left\lceil
 \frac{4\kappa_{\psi,j}(x)}{\varepsilon_1}
  \left(\frac{2\kappa_{\psi,j}(x)}{\varepsilon_1}+\frac{1}{3}\right)
  \log\left(\frac{n+1}{t}\right)\right \rceil\right\}\
\eeqn
and
\beqn{Hk}
|\calH_k| \geq\min\left \{ N, \left \lceil
 \frac{4\kappa_{\psi,j}(x)}{\varepsilon_2}
  \left(\frac{2\kappa_{\psi,j}(x)}{\varepsilon_2}+\frac{1}{3}\right)
  \log\left(\frac{2n}{t}\right)\right \rceil\right \}.
\eeqn
The adaptive nature of these sample sizes is apparent in formulae
\eqref{Dkl}--\eqref{Hk}, because they depend on $x$ and $\varepsilon_j$,
which are themselves dynamically updated in the course of the AR$p$DA
algorithm. Depending on the size of $N$, it may clearly be necessary to
consider the whole set $\{1,\dots,N\}$ for small values of $
\{\varepsilon_j\}_{j=0}^2$.  If the cost of evaluating functions
$\psi_i$, for $1\le i\le N$, is comparable for all $i$, the cost of
evaluating $\barf_k (x_k,\omega_k)$ amounts to the fraction 
$|{\calD_{k,1}}|/N$ of the effort for computing the exact value
$f(x_k)$. Analogous considerations hold for the objective function's
derivatives.

The implementation of rules \eqref{Dkl}-\eqref{Hk} requires the knowledge of
the size of the functions $\psi_i$'s and their first and second order
derivatives. If only global information is available, the dependence on $x$
may obviously be avoided by choosing a uniform upper bound $\kappa_{\psi,j}$
for all $x \in \calF$, at the cost of a lesser adaptivity.  Similar bounds on
the sample size used to approximate gradients and Hessians up to a prescribed
probability have been derived and used in \cite{RoosMaho18} where it has also
been observed that there are problems where estimations of the needed uniform
upper bounds can be obtained. In particular, let $\{(a_i, b_i)\}_{i=1}^N$
denote the pairs forming a data set with $a_i \in \Re^n$ being the vector
containing the features of the $i$-th example and $b_i$ being its label.  In
\cite{RoosMaho18} authors considered the minimization of objective function
$(1/N) \sum_{i=1}^N (\Phi(a_i^Tx)-b_ia_i^Tx)$ over a sparsity inducing
constraint set, e.g., ${\cal F} = \{x \in \Re^n \mid \|x\|_1\le 1\}$, for
cumulant generating functions $\Phi$ of different forms, and explicitly
provided the uniform bound $\kappa_{\psi,1}$.  Taking into account that $x$
belongs to the set ${\cal F} $, uniform bounds for the objective function and
the Hessian norm can also be derived.

Uniform bounds are available also in the unconstrained setting for binary
classification problems modelled by the sigmoid function and least-squares
loss, i.e. problems of the form \req{finite-sum} with ${\cal F}\equiv \Re^n$
and
\beqn{sigmoid-least-squares}
\psi_i(x)=\left(b_i-\frac{1}{1+e^{-a_i^Tx}} \right)^2,\quad i=1\ldots,N.
\eeqn
Let $v_i(x)=(1+e^{-a_i^Tx})^{-1}$ and note that $b_i\in\{0,1\}$ and
$v_i(x)\in (0,1)$ for any $x\in \Re^n$. Then, $|\psi_i(x)|\le 1$, for any
$x\in \Re^n$.  Moreover, uniform upper bounds $\kappa_{\psi,j}$ for
$\nabla_x^j \psi_i(x)$, $j=1,2$ can be easily derived and they are reported
Table \ref{Table_phi_classification} along with the expression of the first
and second order derivatives of $\psi_i(x)$. The computation of these bounds
requires a pre-processing phase as the norm of the features vectors
$\{a_i\}_{i=1}^N$ of the data sets are needed.

\begin{table}[htb]
\begin{center}
\begin{tabular}{l|c|c}
                      & Derivatives & $\kappa_{\psi,j}$ \\ \hline
$\nabla_x \psi_i(x)$   &  $-(b_i-v_i(x))(1-v_i(x))v_i(x) a_i$ & $\sfrac{1}{5} \|a_i\| $\\
$\nabla^2_x \psi_i(x)$ &  $v_i(x)(1-v_i(x))(3v_i(x)^2-2v_i(x)(1+b_i)+b_i)a_ia_i^T$
                       &  $\sfrac{1}{10} \|a_i\|^2$ \\ \hline
\end{tabular}
\caption{First and second order derivatives  of
    \eqref{sigmoid-least-squares} and corresponding uniform
    bounds}\label{Table_phi_classification} 
\end{center}
\end{table}
   
Finally, whenever $N$ is large enough to ensure that \eqref{Dkl}--\eqref{Hk}
do not require the full sample, the size of the sample used to obtain a single
approximate objective function value is $O(\varepsilon_0^{-2})$.  Analogously,
gradient and Hessian values are approximated by averaging over samples of size
$O(\varepsilon_1^{-2})$ and $O(\varepsilon_2^{-2})$, respectively. In Step 3
of the AR$1$DA algorithm, the choice $\varepsilon_0\in \big [\gamma_\epsilon
  \omega_k\|\overline{\nabla_x^j f}(x)\|^2/\sigma_k,
  \omega_k\|\overline{\nabla_x^j f}(x)\|^2/\sigma_k\big]$ is required to
ensure that \eqref{Df+-DT-first}-\eqref{Df-DT-first} are satisfied. With this
choice, iteration $k$ of the AR$1$DA algorithm requires
$O(\|\overline{\nabla_x^j f}(x)\|^{-4})$ $\psi_i$-evaluations
($O(\epsilon^{-4})$ $\psi_i$-evaluations in the worst case).  Similarly,
$\varepsilon_0=O(\omega_k \| \barDT_2^f(x_k,s_k,\omega_k)\|)$ is needed at
iteration $k$ of the AR$2$DA algorithm. As a consequence, and if the algorithm
does not terminate at iteration $k+1$, it follows from Lemma \ref{Dm-lemma}
and \ref{longs-j-lemma} that $O(\epsilon^{-(3/(3-q))^2})$ $\psi_i$-evaluations
may be required in the worst case.  Finally, using Lemma \ref{one-it-compl-l}
and \eqref{wrong-xi-2}, we claim that each iteration of the AR$1$DA and
AR$2$DA algorithms requires at most $O((1+\nu_{\max}(\epsilon))\epsilon^{-2})
$ evaluations of component gradients and component Hessians, where
$\nu_{\max}(\epsilon)$ has been defined in \req{nuk-bound}.  These bounds turn
out to be better or the same as those derived in
\cite{BlanCartMeniSche16},\cite{CartSche17},\cite{YaoXuRoosMaho18}. Although
they may appear discouraging, it should be kept in mind that they are valid
only if $N$ is truly large compared with $1/\epsilon$ (for instance, it has to
exceed $O(\epsilon^{-4})$ to allow for approximate functions in the AR1DA
Algorithm). In other words, the sampling schemes \eqref{Dkl}--\eqref{Hk} are
most relevant when $1/\epsilon$ remains modest compared with $N$.

We conclude by emphasizing that the per-iteration failure probability $t$ given in
\req{t} is not too demanding in what concerns the sample size, because it only
occurs in the logarithm term of \req{sample-size}. The same is true of the
impact of the value of the unknown constants hidden in the $O(\cdot)$ notation
in \req{t}.

\numsection{Conclusion and perspectives}\label{conclusion-s}

We have provided a general regularization algorithm using inexact function and
derivatives' values, featuring a flexible adaptive mechanism for specifying
the amount of inexactness acceptable at each iteration.  This algorithm,
inspired by the unifying framework proposed in \cite{CartGoulToin18b}, is
applicable to unconstrained and inexpensively-constrained nonconvex
optimization problems, and provides optimal iteration complexity for arbitrary
degree of available derivatives, arbitrary order of optimality and the full
range of smoothness assumptions on the objective function highest derivative.
We have also specialized this algorithm to the cases of first- and
second-order methods, exhibiting simple and numerically realistic methods. We
have finally provided a probabilistic version of the complexity analysis and
derived associated lower bounds on sample size in the context of subsampling
methods.

There are of course many ways in which the proposed algorithm might be
improved. For instance, the central calculation of relatively accurate Taylor
increments may possibly be made more efficient by updating the absolute
accuracies for different degrees separately. Further techniques to avoid
unnecessary derivative computations (without affecting the optimal complexity)
could also be investigated.

The framework proposed in this paper also offers obvious avenues for
specializations to specific contexts, among which we outline two.  The first
is that of algorithms using stochastic approximations of function values and
derivatives.  The technique presented here derives probabilistic conditions
under which properties of the deterministic algorithms are preserved.  It does
not provide an algorithm which is robust against failures to satisfy the
adaptive accuracy requirements. This is in contrast with the interesting
analysis of unconstrained first-order methods of \cite{PaquSche18} and
\cite{BlanCartMeniSche16}. Combining the generality of our approach with the
robustness of the proposal in these latter papers is thus desirable. The
second interesting avenue is the application of the new results to
multi-precision optimization in the context of very high performance
computing. In this context, it is of paramount importance to limit energy
dissipation in the course of an accurate calculation, and this may be obtained
by varying the accuracy of the most crucially expensive of its parts (see
\cite{GratSimoToin18} for unconstrained quadratic optimization). The
discussion above again provides guidance at what level of arithmetic accuracy
is needed to achieve overall performance while maintaining optimal
complexity. Both these topics are the object of ongoing research and will be
reported on at a later stage.

{\footnotesize
\section*{\footnotesize Acknowledgments}
%\vspace*{-3mm}

INdAM-GNCS partially supported the first and third authors under Progetti di
Ricerca 2018.  The second author was partially supported by INdAM through a
GNCS grant. The last author gratefully acknowledges the support and friendly
environment provided by the Department of Industrial Engineering at the
Universit\`{a} degli Studi, Florence (Italy) during his visit in the fall of
2018. The authors are also indebted to two careful referees, whose comments
and perceptive questions have resulted in a significant improvement of the
manuscript.

}
\end{document}